\documentclass[journal]{IEEEtran}
\ifCLASSINFOpdf
\else
\fi

\IEEEoverridecommandlockouts
\usepackage{array}
\usepackage{xfrac}
\usepackage{amsfonts}
\usepackage{mathrsfs}
\usepackage{epsfig,latexsym}
\usepackage{amsmath}
\usepackage{amsfonts,amssymb,color}
\usepackage{amssymb,latexsym,amsfonts,amsmath,bbm}
\usepackage{bm}
\usepackage{float}
\usepackage{flushend}
\usepackage{graphicx} 
\usepackage{colortbl}
\definecolor{pinegreen}{rgb}{0.29,0.52,0.42}
\usepackage{flushend}
\allowdisplaybreaks[4]

\newtheorem{thm}{Theorem}
\newtheorem{lem}{Lemma}

\newtheorem{assum}{Assumption}
\newtheorem{rem}{Remark}

\begin{document}
%
\title{Distributed Observer Design over Directed Switching Topologies}
%
%
%

\author{Haotian~Xu,
	    Shuai~Liu,~\IEEEmembership{Member,~IEEE},
        Bohui~Wang,~\IEEEmembership{Senior Member,~IEEE},
        Jingcheng~Wang
\thanks{This work is supported by National Natural Science Foundation of China (No.61821004, 62133008, 61903290). The corresponding author: Shuai Liu.}
\thanks{H. Xu and S. Liu are with School of Control Science and Engineering, Shandong University, Jinan 250061, China~(email: xuhaotian@sdu.edu.cn, liushuai@sdu.edu.cn)}	
\thanks{B. Wang is with the School of Aerospace Science and Technology, Xidian University, Xi'an 710071, China (e-mail: wang31aa@126.com; wshchen@126.com).}
\thanks{J. Wang is with Department of Automation, Shanghai Jiao Tong University, Key Laboratory of System Control and Information Processing, Ministry of Education of China, Shanghai, 200240; He is also with Shanghai Engineering Research Center of Intelligent Control and Management, Shanghai 200240, China~(e-mail: jcwang@sjtu.edu.cn)}}

%
%

\markboth{Submitted to IEEE Transactions on Automatic Control}%
{Haotian Xu \MakeLowercase{\textit{et al.}}: Distributed Observer Design}
%



\maketitle

\begin{abstract}
The distributed observer design problem holds significant importance in cases in which the output information of a system is decentralized across different subsystems. Each subsystem has a local observer and access to one part of the measurement outputs and information exchanged through communication networks. This paper focuses on the design of distributed observer with jointly connected directed switching networks. The problem presents challenges due to passive switching modes and the open-loop unboundedness that results from local observability. To overcome these challenges, we develop a network transformation mapping method whereby each local observer can classify itself into an independent subgraph based on independent judgment. Next, an observable decomposition and reorganization method is developed for the digraph case to ensure that each subgraph possesses independent dynamic properties. Asymptotic omniscience is then proven using a developed recursive proof method. This paper includes many previous results as special cases, because most are only suitable for undirected switching topologies or fast-switching cases. An adaptive coupling gain design is proposed to simplify the calculation and verification of conditions that guarantee asymptotic omniscience. Finally, simulation results with the power system show the validity of the developed theory.
\end{abstract}


%
\IEEEpeerreviewmaketitle

\section{Introduction}
Distributed observer has received extensive attention in recent years  \cite{Xu2021IJRNC,9461598,8093658,2019Completely,battilotti2019distributed} owing to its great potential for improving the performance of distributed control, such as the distributed cooperative control of battery energy storage systems in DC microgrids  \cite{Meng2021Distributed}, distributed controller designed for flexible structure \cite{Zhang2018Distributed} and distributed Nash equilibrium game for nonlinear systems \cite{Huang2020Distributed}.

Research on distributed observer is divided into two categories. The first focuses on estimating the states of an external system and the second on estimating the states of an interconnected system in a distributed way \cite{Xu2020AnSMC}. Since the first kind of distributed observer is the dual version of leader-following systems \cite{Zhang2011Optimal,2014The,Huang2017The}, it is regarded as a kind of multi-agent system. The distributed observer mentioned later refers to the second one. 

Consider the linear system $\dot{\chi}=A\chi,~y=C\chi$ in which $A,C$ have compatible dimensions and $\chi\in\mathbb{R}^n, y=[y_1^T,\ldots,y_N^T]^T\in\mathbb{R}^m$ stand for states and measurement outputs respectively. Distributed observer is a cooperative observer network composed of multiple local observers (also termed ``agents"). Each local observer is assumed to have access to one part of outputs $y_i$ and the information exchanged via a communication network. Let $\hat{\chi}_i$ be the estimates of $\chi$ generated by the $i$th agent. Then, distributed observer aims to achieve omniscience asymptotically, i.e., $\lim_{t\to\infty}\|\hat{\chi}_i(t)-\chi(t)\|=0$ for all $i=1,\ldots,N$. 

Up to now, distributed observer has been widely studied, such as the earliest distributed Luenberger observer \cite{Han2017A}, distributed minimum order observer \cite{Han2018Towards}, distributed observer with event-triggering mechanism \cite{9461598}, communication delay \cite{Liu2018Design}, network attack \cite{deghat2019detection}, adaptive coupling gain \cite{2019Completely}, and so on. Distributed-observer-based control approaches have been studied in \cite{2017CooperativeTCS,Liu2018Cooperative,Xu2020Distributed}, and some claim that the distributed control law can achieve a control performance similar to the centralized control \cite{Zhang2018Distributed,Huang2020Distributed,Xu2022TIV}. It also provides profound application value for research on distributed observer.

In the past few years, the problem of the distributed observer under switching topologies has attracted much attention because failure or occasional interruption of the communication network is inevitable \cite{YANG2021539,deghat2019detection,mitra2019byzantine}. One of the significant challenges described in \cite{Xu2021TCYB} is that local state estimates of several different local observers cannot achieve consensus under the jointly connected topologies (which are disconnected at any time), even if they belong to the same connected branch. It is significantly different from the multi-agent problem \cite{Cheng2020Seeking}. {\color{blue}Also, the leader-following consensus of multi-agent systems under real-time disconnected switching topologies relies on a marginally stable open-loop leader system (refer to \cite{HE2021110021}). However, its approaches are unsuitable for distributed observer because the non-consistency terms in distributed observer are generally unstable.} Consequently, most previous results avoid real-time disconnected switching topologies and turn to the case in which augmented topologies contain strongly connected networks \cite{Wanglili2020A} or switch in fast-switching mode \cite{Wen2014Consensus,2017CooperativeTCS,2017Cooperative}.

{\color{blue}Recently, some studies have investigated the real-time disconnected switching topologies problem in distributed observer design without resorting to a fast-switching approach. \cite{zhang2023distributed} achieves distributed observer with a jointly connected network, but it assumes that the system matrix $A$ is marginally stable. However, unlike the multi-agent problem, we usually do not assume that the estimated system is open-loop marginally stable in the distributed observer problem, because the stability of the observer-based control system needs to be proved---and it is thus impossible to know in advance whether the observed system is marginally stable. Two studies, \cite{Xu2020IFAC} and \cite{YANG2023110690}, have implemented the distributed observer design for a class of unstable $A$ (relying on some constraint conditions). However, their stability proof depends on the approach developed by \cite{Wei2010Leader}, and this approach is only applicable to the case of undirected graphs and cannot be extended to directed graphs.}

Therefore, when the system matrix is unstable and the communication topology is jointly connected, how to design and implement a distributed observer with asymptotic omniscience is unclear. Studies by \cite{Xu2021TCYB} have shed light on this issue by requiring that the output matrix $C$ be known by all agents and addressing a states decomposition and reorganization method. With this approach, asymptotic omniscience under jointly connected networks can be achieved without any system constraints. 

{\color{blue}The fly in the ointment is that \cite{Xu2021TCYB} only focuses on undirected communication networks. Transformation from an undirected graph to a directed graph will pose significant challenges in the case in which networks are not strongly connected at any time. This is significantly different from the fixed topology scenario, such as from the undirected graph case \cite{Kim2016Distributed} to the directed graph case \cite{Han2017A}. The reason is that unconnected undirected graphs always contain connected branches, yet directed graphs do not in general include a counterpart.} Specifically, \cite{Xu2021TCYB} uses the state reorganization and decomposition method in the undirected graph to construct a group of novel distributed observers that can realize omniscience asymptotically in each connected branch. However, there is no counterpart in a directed graph the corresponds to the connected branch of an undirected graph. In most cases, there is not even a separate subgraph in a directed graph. Even if separate subgraphs exist, their nodes cannot share information sufficiently because they are not always strongly connected, which prevent use of the state reorganization and decomposition method. This will undoubtedly force researchers to develop new approaches and devise state decomposition and reorganization methods in directed graph cases. 

This paper overcomes these difficulties in three steps. First, we consider a network transformation mapping such that the transformed topology contains an independent subgraph. Moreover, network transformation mapping can be directly calculated by each local observer using the topology information of the current network. Second, the states decomposition and reorganization method for the digraph case is developed with the assistance of network transformation mapping. This allows the dynamics of each local observer to be decomposed into several different observable subspaces. Then, different observable subsystems of each local observer employ different transformed networks. As a result, the system after state reorganization has relatively independent dynamic performance and there is a lower triangular interconnection relationship between each group of reorganized systems. Finally, a recursive proof method of state omniscience is proposed. Through the convergence of error dynamics about the first group of reorganized states and the lower triangular interconnection relationship, the convergence of all reorganized error dynamics can be proved iteratively. Then, the asymptotic omniscience of the distributed observer is obtained by applying the state decomposition and reorganization method in reverse.

The contributions of this paper are concluded as follows. 1) Asymptotic omniscience is achieved based on the proposed distributed observer under jointly connected switching directed networks. In comparison, methods in the literature usually rely on the undirected graph case \cite{Xu2021TCYB}, fast switching (e.g., \cite{2017CooperativeTCS,2017Cooperative}), and real-time strongly connected networks \cite{Wanglili2020A}. 2) An original proposed network transformation mapping is investigated such that the transformed network contains an independent subgraph with spanning trees. This solves the problem whereby there are generally no independent subgraphs in a directed graph. 3) With the help of network transformation mapping, this paper proposes a reconstructed states decomposition and reorganization method for directed graph such that each independent subgraph can contain relatively independent observer dynamics. Similar methods are only suited for undirected graphs in the literature. 4) With the aid of the recursive proof method, fix and adaptive coupling gains are proposed to guarantee asymptotic omniscience.

The rest of th paper is organized as follows. The problem is formulated in Section \ref{sec2}, and the method for designing distributed observer is presented in Section \ref{sec3}. Section \ref{sec4} demonstrates the stability of the error dynamics with the states decomposition and reorganization method and recursive proof method. Section \ref{sec5} supplements these results in terms of adaptive coupling gain. Two numerical simulations with respect to the power system are employed in Section \ref{sec6} to verify the stated theories, and Section \ref{sec7} concludes the paper.



\section{Problem formulation}\label{sec2}
Notation used throughout the paper are as follows. 

\textbf{(N1)} Denote $x\in\mathbb{R}^n$ and $A\in\mathbb{R}^{n\times n}$ the $n$ dimension column vector and $n\times n$ dimension matrix, respectively. $x^T, A^T$ stand for their transpositions. In particular, $I_n$ is the $n$ dimension identity matrix. $\|\cdot\|$ denotes the vector or the matrix norm induced by $2$-norm. Given arbitrary elements $a_i, i=1,\ldots,m$, $col\{a_1,\ldots,a_m\}\in\mathbb{R}^m$ represents the stack of $a_i, i=1,\ldots,m$ in a column and $diag\{a_1,\ldots,a_m\}\in\mathbb{R}^{m\times m}$ is the block diagonal matrix with all $a_i, i=1,\ldots,m$ on its diagonal. Denote $\bar{\lambda}(\cdot)$ and $\underline{\lambda}(\cdot)$ the maximum eigenvalue and minimum eigenvalue of a symmetric matrix, respectively. $sym\{A\}$ stands for $A+A^T$ for $A\in\mathbb{R}^{n\times n}$. Operators $\wedge$ and $\vee$ are defined as $a\wedge b=\min\{a,b\}$ and $a\vee b=\max\{a,b\}$ for $a,b\in\mathbb{R}$.

\textbf{(N2)} Set $\mathcal{G}=\{\mathcal{V},\mathcal{E},\mathcal{A}\}$ to be a directed graph and $\mathcal{V}$ and $\mathcal{E}$ the set of its nodes and arcs respectively, and denote by $\mathcal{A}$ its adjacency matrix. The element $a_{ij}=1$ of $\mathcal{A}$ indicates an arc pointing from node $j$ to $i$, and $a_{ij}=0$ otherwise. A directed path from $i_1$ to $i_k$ composed of a set of nodes $\{i_1,\ldots,i_k\}$ is defined by $a_{i_{p+1}i_p}=1$ for all $p\in\{1,\ldots,k\}$. We say $i_k$ is able to obtain information from $i_1$ by $l$ hops if the path between $i_k$ and $i_1$ includes $l$ arcs. An input-degree matrix $D=diag\{d_1,\ldots,d_N\}$ is defined by $d_i=\sum_{j=1}^Na_{ij}$, and the matrix $\mathcal{L}=D-\mathcal{A}$ is the graph Laplacian of $\mathcal{G}$. The directed graph is a strongly connected graph if there is a directed path between any pair of nodes belonging to $\mathcal{V}$. {\color{blue}The network contains a spanning tree if there is a node $k$ and the directed path pointing from $k$ to $i$ exists $\forall i\in\mathcal{V}$ and $i\neq k$, where $k$ is said to be a root of the spanning tree. }

Consider a linear system,
\begin{align}
	&\dot{\chi}=A\chi,\\
	&y=C\chi,
\end{align}
where $\chi\in\mathbb{R}^n$ and $y\in\mathbb{R}^p$ are the states and measurement outputs respectively; system matrix $A$ and output matrix $C$ are with compatible dimension. Assume that the pair $(C,A)$ is observable and further assume that each of $N$ agents only has access to one part of $y$. Also, they are capable of reconstructing complete output information $y$ by collaborating. The original output vector and output matrix can be rewritten as $y=col\{y_1,\ldots,y_N\}$ and $C=col\{C_1,\ldots,C_N\}$, respectively. $y_i=C_i\chi\in\mathbb{R}^{p_i}$ stands for the information obtained by the $i$th agent (or the ``$i$th node" in the communication network), where $p=\sum_{i=1}^Np_i$. The pair $(C_i,A)$ does not need to be observable for all $i=1,\ldots,N$. In this paper, the output matrix $C$ is assumed to be known for all agents. Hence, {\color{red} there exists a nonsingular matrix $T$ that leads to a lower triangular form}:
\begin{align}
	&\dot{x}_{0,i}=A_{io}x_{0,i}+\sum_{l=1}^{i-1}\Upsilon_{il}x_{0,l},\label{decom1}\\
	&y=\Sigma x=col\{\Sigma_1,\ldots,\Sigma_N\}x,\label{decom2}
\end{align}
where $v_i$ is the observability index of the $i$th observable subsystem; $x_{0,i}\in\mathbb{R}^{v_i}$ is the element of $x=col\{x_{0,1},\ldots,x_{0,N}\}$ with $x=T\chi$; $A_{io}\in\mathbb{R}^{v_i\times v_i}$, $C_{io}\in\mathbb{R}^{p_i\times v_i}$, $\Upsilon_{il}\in\mathbb{R}^{v_l\times v_i}$; and
\begin{align}
	\Sigma_i=\begin{bmatrix}0_{p_i\times\bar{v}_{i-1}},C_{io},0_{p_i\times(n-\bar{v}_i)}\end{bmatrix}
\end{align}
are matrices transformed from $A$ and $C$ with $\bar{v}_i=\sum_{j=1}^iv_i$. {\color{red}Moreover, $\Upsilon_{il}$ stands for the coupling matrix between $x_{0,l}$ and $x_{0,i}$; $v_i$ and $v_l$ are the dimensions of $A_{io}$ and $A_{lo}$, respectively.} Note that the pair $(C_{io},A_{io})$ is observable according to the property of observability decomposition. {\color{red}To further help the reader understand the relationship between $A$ and $A_{io}$, $\Upsilon_{il}$, we specifically show that
\begin{align}
	TAT^{-1}=\begin{bmatrix}
		A_{1o}&&&\\ \Upsilon_{21}&A_{2o}&&\\\vdots&\vdots&\ddots&\\\Upsilon_{N1}&\Upsilon_{N2}&\cdots&A_{No}
	\end{bmatrix}.
\end{align}
}

The overall goal of this paper is to design a distributed observer for the underlying system. The $i$th agent owns the outputs of the $i$th observable subsystem and establishes the $i$th local observer, such that its local state estimation $\hat{\chi}_i$ satisfies $\lim_{t\to\infty}\|\hat{\chi}_i-\chi\|=0$, or $\lim_{t\to\infty}\|\hat{x}_i-x\|=0$, where $\hat{x}_i=T\hat{\chi}_i$. As stated in Section 1, this property is the so-called \textit{asymptotic omniscience} defined by \cite{8093658} and the references therein. 

The communication network among $N$ agents will be described as follows, and each agent is regarded as a node in the network. Given a piecewise function $\sigma(t)$ that maps time $t$ to an index set $\mathcal{P}=\{1,\ldots,\varpi\}$. The network that appears in time interval $[0,+\infty)$ is indexed by $\sigma(t)$. For example, a graph at time $t$ can be characterized as $\mathcal{G}^{\sigma}=\{\mathcal{V}^\sigma,\mathcal{E}^\sigma\}$, where $\mathcal{V}^\sigma$ and $\mathcal{E}^\sigma$ are its node set and arc set, respectively. Let $\mathcal{A}^\sigma=[\alpha_{ij}^\sigma]_{i,j=1}^N$ be the adjacency matrix of $\mathcal{G}^\sigma$. We now state the following two assumptions.

\begin{assum}\label{assume3}
	Consider a {\color{red}time series $t_0,~t_1,\cdots$, $t_k,\cdots$ }subject to $t_{k+1}-t_k= \mathcal{T}$ for $k=0,1,\cdots$, and a subsequence in $[t_k,t_{k+1})$ defined as $t_k^0,\cdots,t_k^{\varpi}$ with $t_k^{j+1}-t_k^j=\tau=\mathcal{T}/\varpi, j=1,\cdots,\ell_k-1$, where $\mathcal{T}>0$ is a given constant, {\color{red}$t_0$ is the initial time of the system and $t_j$ is the time moment with respect to the switching.} Then, the interconnection topology is supposed to remain unchanged in each $[t_k^j,t_k^{j+1})$, and it is further supposed that $\mathcal{G}^{\sigma(t)}=\mathcal{G}^j$ if $t\in[t_k^j,t_k^{j+1})$ for arbitrary $k=0,1,\cdots$.
\end{assum}

\begin{assum}\label{assume2}
	We assume that the union $\bigcup_{\sigma\in\mathcal{P}}\mathcal{G}^\sigma$ of all $\mathcal{G}^\sigma,~\sigma\in\mathcal{P}$ is strongly connected. Besides, $\forall i,j$, there exists at least one graph belonging to $\mathcal{P}$ which contains a path from $i$ to $j$.
\end{assum}

\begin{rem}
Both Assumption \ref{assume3} and Assumption \ref{assume2} are mild. Assumption \ref{assume3} differs from the traditional assumption about jointly connected switching topologies, which it requires a fixed dwell time. This is in line with some practical situations. For example, each node switches communication signals according to a certain communication protocol if communication resources are limited. Assumption \ref{assume2} has one more directed path condition than the traditional assumption. The directed path condition holds when network resources are limited and some links in the strongly connected graph need to be interrupted intermittently. {\color{blue}Furthermore, we will give some illustration concerning the assumption of ``knowing the output matrix $C$ for all agents". Indeed, many interconnected systems in which the output functions of each subsystem are the same, such as the power system consisting of several power-generation areas, the droop control systems of microgrids, and the coupled Van der Pol system. In these systems, as long as the agent knows the form of its own output function, this is equivalent to knowing the form of other subsystems' output functions. Therefore, from this perspective, assuming that matrix $C$ is known to all agents does not compromise the distributed characteristics of the distributed observer.}
\end{rem}

\section{\color{pinegreen}Design of distributed observer}\label{sec3}

{\color{pinegreen}
In this section, we first present dynamic equations for the distributed observer in Subsection \ref{sec3.1}. Network transformation mapping will be introduced in Subsection \ref{sec3.2}, which is directly related to the acquisition of some key parameters in distributed observer. Also, network transformation mapping is one of the key technologies this paper can use to achieve distributed observer for open-loop unstable systems under joint connected switching topology conditions.}

\subsection{Dynamics of the local observer}\label{sec3.1}

Dynamics of the $i$th local observer take the form of
\begin{align}
	\dot{\hat{x}}_{i,i}=&A_{io}\hat{x}_{i,i}+H_{io}\left(y_i-C_{io}\hat{x}_{i,i}\right)\notag\\
	&+\mathbbm{r}_{i,i}\sum_{j=1}^N\alpha_{ij}^{\sigma}(i)\left(\hat{x}_{j,i}-\hat{x}_{i,i}\right)+\sum_{l=1}^{i-1}\Upsilon_{il}\hat{x}_{i,l},\label{do1}\\
	\dot{\hat{x}}_{i,k}=&A_{ko}\hat{x}_{i,k}+\mathbbm{r}_{i,k} \sum_{j=1}^N\alpha_{ij}^{\sigma}(k)\left(\hat{x}_{j,k}-\hat{x}_{i,k}\right)\notag\\
	&+\sum_{l=1}^{k-1}\Upsilon_{kl}\hat{x}_{i,l},~~k=1,\ldots,N,~k\neq i.\label{do2}
\end{align}
Herein, coupling gains $\mathbbm{r}_{i,k}$ are pending parameters. $H_{io}$ is the observer gain chosen so that $A_{io}-H_{io}C_{io}$ is a Hurwitz matrix. $\alpha_{ij}^{\sigma}(k)$, which is supposed to be highlighted, is the image of the network transformation. In the next section, we will introduce the construction method of $\alpha_{ij}^{\sigma}(k)$ and its function.

First, we state two lemmas (see the proof of Lemma \ref{cyb} in the Appendix) that will be used in this paper. 
\begin{lem}\label{cyb}
	Consider a linear time-variant system
	\begin{equation*}\label{eq-lem1}
		\dot{x}(t)=A(t)x(t)+M\xi(t),
	\end{equation*}
	where $x\in\mathbb{R}^n$, $A(t),~M\in\mathbb{R}^{n\times n}$, and $\xi(t)\in\mathbb{R}^n$ satisfies $\lim_{t\to\infty}\xi(t)=0$, then $\lim_{t\to\infty}x(t)=0$ if the solution $\mathbbm{x}(t)$ of $\dot{x}(t)=A(t)x(t)$ satisfies $\mathbbm{x}(t)\leq a_1e^{-a_2(t-t_0)}$, where $a_1,a_2$ are positive constants and $t_0>0$ is the initial time.
\end{lem}
\begin{lem}[\cite{Hong2008Distributed}]\label{spanning}
	Consider a graph $\mathcal{G}=\{\mathcal{V},\mathcal{E},\mathcal{A}\}$. Assume $\mathcal{G}$ contains a spanning tree rooted at $i_r\notin\mathcal{V}$. Set $b_i=1$ if there is a path from $i_r$ to $i\in\mathcal{V}$. Then, all eigenvalues of matrix $\mathcal{H}=\mathcal{L}+\mathcal{B}$ locate at the open right half of the plane, where $\mathcal{B}\triangleq diag\{b_i,~i\in\mathcal{V}\}$.   
\end{lem}

\begin{figure*}[!t]
	\centering
	\includegraphics[width=16cm]{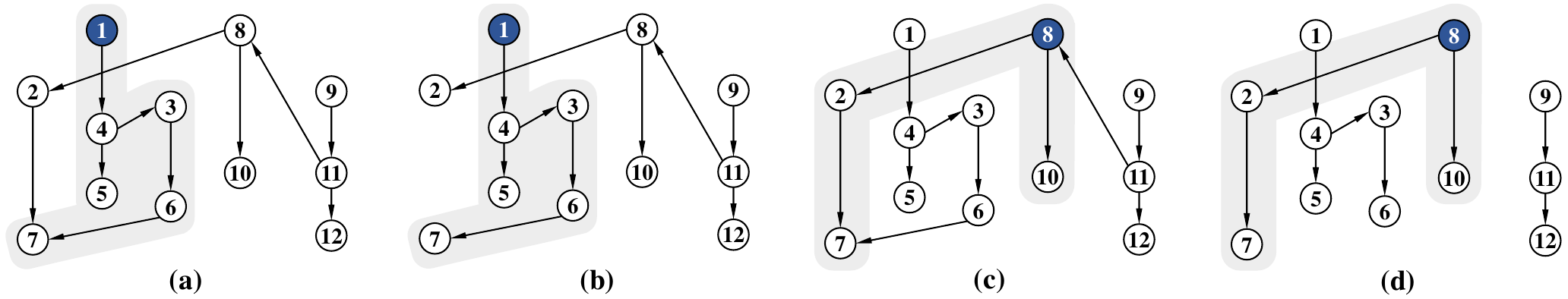}\\
	\caption{Illustration of the mapping on an adjacency matrix}\label{zhuanhuantu}
\end{figure*}

\subsection{Network transformation mapping}\label{sec3.2}

Generally speaking, there is no strongly connected branch in directed switching topologies, which means that the state decomposition and reorganization method based on connected branches proposed in \cite{Xu2021TCYB} is not suitable for directed graphs. Nevertheless, a directed graph always has a subgraph composed of some nodes with a tree structure. Although this tree structure does not have good network properties like that of spanning tree owing to disturbance from nodes outside the tree, it offers a new sight to developing state decomposition and reorganization methods for directed graph cases. However, because the tree structure lack the excellent properties of a spanning tree, it cannot guarantee the asymptotic omniscience of the distributed observer. Therefore, this section is devoted to proposing a network transformation mapping to construct a subgraph that contains a spanning tree. 

To this end, denote $l_{ik}^{\sigma}$ the element located at row $i$ and column $k$ of the matrix  $\mathcal{A}^{\sigma}_{r}=I_N+\mathcal{A}^{\sigma}+\left(\mathcal{A}^{\sigma}\right)^2+\cdots+\left(\mathcal{A}^{\sigma}\right)^{N-1}$. Then, let $\ell_{ik}^\sigma=l_{ik}^\sigma\wedge 1$ be the elements of $\mathcal{A}^{\sigma}_{r}$'s structure matrix, i.e, $\ell_{ik}^\sigma=1$ if $l_{ik}^\sigma>0$ and $\ell_{ik}^\sigma=0$ otherwise. Construct a node set with respect to a node $k$ as $\mathcal{V}_k^\sigma=\{i\in\mathcal{V}:~\ell_{ik}^\sigma=1\}\cup\{k\}$. Hence, we can state the following.

\begin{lem}\label{renet}
	Network $\mathcal{G}^\sigma$ is transformed to $\mathcal{G}^\sigma(k)$ with respect to an arbitrary given node $k$ if the elements of its adjacency matrix $\alpha_{ij}^{\sigma}$ are transformed by
	\begin{align}
		&\alpha_{ij}^{\sigma}(k)=\delta_{ij}^{\sigma}\alpha_{ij}^{\sigma},\label{mapping1}\\
		&\delta_{ij}^{\sigma}=\left(\left(1-\ell_{ik}^\sigma\right)\wedge\left(1-\ell_{jk}^\sigma\right)\right)\vee\left(\ell_{ik}^\sigma\wedge\ell_{jk}^\sigma\right).\label{mapping2}
	\end{align}
	Then, $\mathcal{G}^\sigma(k)$ is the transformed network with adjacency matrix $[\alpha_{ij}^{\sigma}(k)]_{i,j=1}^N$, and contains an independent subgraph $\mathcal{G}^{\sigma}_{k}\triangleq\{\mathcal{V}^{\sigma}_{k},\mathcal{E}^{\sigma}_{k}\}$. In addition, $\mathcal{G}^\sigma(k)$ has a spanning tree rooted at $k$. Furthermore, all the nodes belonging to $\mathcal{V}_{k}^\sigma$ have no access to the information from $\mathcal{V}\backslash\mathcal{V}_{k}^\sigma$.
\end{lem}
\begin{IEEEproof}
	Based on the knowledge about complex networks, node $i$ obtains the information of node $j$ at time $t$ by $1$ hop if and only if $\alpha_{ij}^{\sigma}=1$. Moreover, node $i$ can get node $j$'s information by $2$ hops at $t$ if and only if $\sum_{l=1}^N\alpha_{il}^{\sigma}\alpha_{lj}^{\sigma}>0$. Hence, matrix $\mathcal{A}^{\sigma}\cdot\mathcal{A}^{\sigma}$ denotes the reachable matrix within $2$ hops. Similarly, $\left(\mathcal{A}^{\sigma}\right)^l$ denotes the reachable matrix with $l$ hops. 
	
	As a result, the element $\ell_{ij}^\sigma=1$ denotes that node $i$ obtains node $j$'s information within $N-1$ hops. Then, we can construct $\mathcal{G}^{\sigma}_{k}$ by defining $\mathcal{V}^{\sigma}_{k}=\{k\}\cup\{l:~\ell_{lk}^{\sigma}\neq 0\}$ and $\mathcal{E}^{\sigma}_{k}=\{(i,j)\in\mathcal{E}^{\sigma}:~\alpha_{ij}^{\sigma}(k)=1\}$. Therefore, there exists a directed path from $k$ to $i$ for $\forall i\in\mathcal{V}_k^\sigma$.
	
	To show that $\mathcal{G}^{\sigma}_{k}$ contains a spanning tree, we need to discuss the elements in $\mathcal{E}_k^\sigma$ in the following three cases. 1) The case in which $i,j\in\mathcal{V}^{\sigma}_{K}$. The aforementioned proof leads thereby to $\ell_{ik}^{\sigma}=1$ and $\ell_{jk}^{\sigma}=1$. Consequently, 
	\begin{align*}
		\alpha_{ij}^\sigma(k)=\delta_{ij}^\sigma\alpha_{ij}^\sigma(0\wedge 0)\vee(1\wedge 1)=\alpha_{ij}^\sigma.
	\end{align*}
	
	2) The case in which $i,j\in\mathcal{V}\backslash\mathcal{V}_k^\sigma$. In light of the definition of $\mathcal{V}^{\sigma}_{k}$, we know there are no arcs from node $k$ pointing to $\mathcal{V}\backslash\mathcal{V}^{\sigma}_{k}$. Hence, $\ell_{ik}^\sigma=\ell_{jk}^\sigma=0$. This means that $\alpha_{ij}(k)=(1\wedge 1)\vee(0\wedge 0)\alpha_{ij}^\sigma=\alpha_{ij}^\sigma$.
	
	3) The case in which $i\in\mathcal{V}_k^\sigma$ and $j\in\mathcal{V}\backslash\mathcal{V}_k^\sigma$. Since $l_{ik}^{\sigma}=1$ and $l_{jk}^{\sigma}=0$, we have $\alpha_{ij}(k)=(0\wedge 1)\vee(1\wedge 0)\alpha_{ij}^\sigma=0$. 
	
	The above proof indicates that the network transformation mapping leaves the arcs between the	pair belonging to $\mathcal{V}_{k}^\sigma$ and $\mathcal{V}\backslash\mathcal{V}_{k}^\sigma$ as they are and removes the arcs between $\mathcal{V}_{k}^\sigma$ and $\mathcal{V}\backslash\mathcal{V}_{k}^\sigma$. In other words, there is no path from $\mathcal{V}\backslash\mathcal{V}_{k}^\sigma$ to $\mathcal{V}_{k}^\sigma$, and also no path from $\mathcal{V}_{k}^\sigma$ to $\mathcal{V}\backslash\mathcal{V}_{k}^\sigma$. Hence, $\mathcal{G}_k$ is an independent subgraph. Moreover, because there exists a directed path from $k$ to $i$ $\forall i\in\mathcal{V}_k^\sigma$, $\mathcal{G}_k$ contains a spanning tree rooted at $k$.    
\end{IEEEproof}

According to the definition of Lemma \ref{renet}, coupling gain can be preliminary designed as
\begin{align}\label{doubler}
	\mathbbm{r}_{i,k}=(1-\ell_{ik}^{\sigma})\gamma+\ell_{ik}^{\sigma}\gamma_{i,k},~k=1,\ldots,N,
\end{align}
where $\gamma$ and $\gamma_{i,j}$ are positive parameters that need to be designed later.

\begin{rem}\label{rem-ntm}
	This remark employs an example to show the function of Lemma \ref{renet}. In Figure \ref{zhuanhuantu}, the subfigures (a) and (c) are original directed networks. First, transform (a) with respect to node $1$ by (\ref{mapping1}), (\ref{mapping2}). Then, (a) can be transformed into subfigure (b), which contains an independent subgraph $\mathcal{G}_1$ with $\mathcal{V}_1=\{1,3,4,5,6,7\}$. The arcs among $\mathcal{V}_1$ are preserved and the arcs that point from other nodes to $\mathcal{V}_1$ are removed. Furthermore, $\mathcal{G}_1$ contains a spanning tree rooted at node $1$. To further verify the effectiveness of Lemma \ref{renet}, we transform network (c) with respect to node $8$ and thus obtain the transformation in (d). A subgraph $\mathcal{G}_8$ with $\mathcal{V}_8=\{2,7,8,10\}$ does not obtain any information from the nodes outside and includes a spanning tree rooted at $8$.
\end{rem}

\begin{rem}
	As explained in the introduction, multiple local observers cannot achieve consistency even if they are in the same connected branch (undirected graph), which is caused by the strong coupling between observable subsystems. Moreover, the network transformation mapping in this section can only obtain directed independent subgraphs, whose network properties are not as good as the connected branches of undirected graphs. Therefore, to analyze the performance of a distributed observer in independent subgraphs, it is also necessary to propose a state decomposition and reorganization method (Section \ref{sec4}) for directed graph cases, in order to weaken the coupling relationship and cause the independent subgraphs to have independent dynamics.
\end{rem}

\section{\color{pinegreen}Asymptotic omniscience of distributed observer}\label{sec4}

{\color{pinegreen}The definition of asymptotic omniscience suggests that achieving asymptotic omniscience for distributed observer is equivalent to the stability of its error dynamics. To demonstrate the stability of the error dynamics of distributed observer, this section first presents a state decomposition and reorganization approach in the case of directed graph, and derives some preliminary conclusions (in Subsection \ref{sec4.1}). Subsequently, based on these preliminary conclusions, we will prove the stability of error dynamics in Subsection \ref{sec4.2}.}

\subsection{State decomposition and reorganization}\label{sec4.1}

In this section, state decomposition and reorganization methods are presented to assist in analyzing the error dynamic performance of distributed observer (\ref{do1}), (\ref{do2}), bearing in mind that the states have been decomposed in (\ref{decom1}), (\ref{decom2}). This section will first state the observer error dynamics based on state decomposition, then reorganize them on the basis of observable subspace.

The error dynamics are governed by 
\begin{align}
	\dot{e}_{i,i}=&\left(A_{io}-H_{io}C_{io}\right)e_{i,i}+\sum_{l=1}^{i-1}\Upsilon_{il}e_{i,l},\label{do3}\\
	\dot{e}_{i,k}=&A_{ko}e_{i,k}+\mathbbm{r}_{i,k} \sum_{j=1}^N\alpha_{ij}^{\sigma}(k)\left(e_{j,k}-e_{i,k}\right)\notag\\
	&+\sum_{l=1}^{k-1}\Upsilon_{kl}e_{k,l},~~k=1,\ldots,N,~k\neq i.\label{do4}
\end{align} 
Note that the transformed network with respect to node $i$ indicates there is no arc pointing to node $i$, so $e_{i,i}$ does not contain a coupling term. Therefore, error dynamics can be reorganized in the compact form by the following rule. Denote $\varepsilon_{i\star}=col\{e_{j,i},~j\in\mathcal{V}_{i}^{\sigma}\backslash\{i\}\}$, $\varepsilon_{i\diamond}=col\{e_{k,i},~k\in\mathcal{V}\backslash\mathcal{V}_{i}^{\sigma}\}$, and $\varepsilon_{ii}=1_{\pi_i(\sigma)}\otimes e_{i,i}$, where $\pi_i(\sigma)=\left|\mathcal{V}_{i}^{\sigma}\right|-1$; and denote $\xi_{i\star}=col\{\xi_{k,i},~k\in\mathcal{V}_{i}^{\sigma}\backslash\{i\}\}$ with $\xi_{k,i}=\sum_{j=1}^N\alpha_{kj}^{\sigma}(i)\left(\hat{x}_{j,i}-\hat{x}_{k,i}\right)$. To move on, we denote the Laplacian matrix of $\mathcal{V}_i\backslash\{i\}$ as $\mathcal{L}_{\pi_i(\sigma)}$. Then, one states the properties of the reorganized states. 

\begin{lem}\label{lem7}
	The stability of $col\left\{\varepsilon_{ii},\varepsilon_{i\star}\right\}$ is equivalent to the stability of system $col\left\{\varepsilon_{ii},\xi_{i\star}\right\}$ with the form
	\begin{align}\label{varepsilonii}
		\dot{\varepsilon}_{ii}=&\left(I_{\pi_i(\sigma)}\otimes \left(A_{io}-H_{io}C_{io}\right)\right)\varepsilon_{ii}\notag\\
		&+\sum_{l=0}^{i-1}\left(I_{\pi_i(\sigma)}\otimes \Upsilon_{il}\right)\left(1_{\pi_i(\sigma)}\otimes e_{i,l}\right),
	\end{align}
	and
	\begin{align}\label{xi2}
		\dot{\xi}_{i\star}=&\left(I_{\pi_i(\sigma)}\otimes A_{io}\right)\xi_{i\star}-\left(\mathcal{H}_{\pi_i(\sigma)}\Gamma_i^{\sigma}\otimes I_{v_i}\right)\xi_{i\star}\notag\\
		&+\left(B_{i\star}^{\sigma}\otimes H_{io}C_{io}\right)\varepsilon_{ii}\notag\\
		&-\sum_{l=0}^{i-1}\left(B_{i\star}^{\sigma}\otimes I_{v_i}\right)\left(I_{\pi_i(\sigma)}\otimes \Upsilon_{il}\right)\left(1_{\pi_i(\sigma)}\otimes e_{i,l}\right)\notag\\
		&+\sum_{l=0}^{i-1}\left(\mathcal{H}_{\pi_i(\sigma)}\otimes I_{v_i}\right)\Upsilon_{i\star,l}\varepsilon_{i\star,l}.
	\end{align}
\end{lem}
\begin{IEEEproof}
	Based on the definition, it is easy to get (\ref{varepsilonii}). Furthermore, according to (\ref{doubler}) that $\ell_{ji}^{\sigma}=1$ when $j\in\mathcal{V}_{i}^{\sigma}\backslash\{i\}$, we have
	\begin{align}
		\dot{\varepsilon}_{i\star}=&\left(I_{\pi_i(\sigma)}\otimes A_{io}\right)\varepsilon_{i\star}-\left(\Gamma_i^{\sigma}\mathcal{H}_{\pi_i(\sigma)}\otimes I_{v_i}\right)\varepsilon_{i\star}\notag\\
		&+\left(\Gamma_i^{\sigma}\otimes I_{v_i}\right)\left(B_{i\star}^{\sigma}\otimes I_{v_i}\right)\varepsilon_{ii}+\sum_{l=1}^{i-1}\Upsilon_{i\star,l}\varepsilon_{i\star,l},
	\end{align}
	where  $\Gamma_i^{\sigma}=diag\left\{\gamma_{j,i},~j\in\mathcal{V}_{i}^{\sigma}\backslash\{i\}\right\}$, and $\mathcal{H}_{\pi_i(\sigma)}=\mathcal{L}_{\pi_i(\sigma)}+B_{i\star}^{\sigma}$ with $B_{i\star}^{\sigma}=diag\{\alpha_{ji}(i),~j\in\mathcal{V}_{i}^{\sigma}\backslash\{i\}\}$. Also, $\Upsilon_{i\star,t}=diag\{\Upsilon_{jt},~j\in\mathcal{V}_{i}^{\sigma}\backslash\{i\}\}$ and $\varepsilon_{i\star,t}=col\{e_{j,t},~j\in\mathcal{V}_{i}^{\sigma}\backslash\{i\}\}$.
	
	From the notation of $\xi_{i\star}$, we know
	\begin{align}
		\xi_{i\star}=\left(\mathcal{H}_{\pi_i(\sigma)}\otimes I_{v_i}\right)\varepsilon_{i\star}-\left(B_{i\star}^{\sigma}\otimes I_{v_i}\right)\varepsilon_{ii}.
	\end{align}
	Obviously, based on the stability of $\varepsilon_{ii}$, $\xi_{i\star}$ converges to zero if and only if $\varepsilon_{ii}$, $\varepsilon_{i\star}$ converges to zero. Then, in the time interval in which $\sigma$ is a constant, we have
	\begin{align}\label{xi}
		\dot{\xi}_{i\star}=&\left(\mathcal{H}_{\pi_i(\sigma)}\otimes I_{v_i}\right)\dot{\varepsilon}_{i\star}-\left(B_{i\star}^{\sigma}\otimes I_{v_i}\right)\dot{\varepsilon}_{ii}\notag\\
		=&\left(\mathcal{H}_{\pi_i(\sigma)}\otimes I_{v_i}\right)\left(I_{\pi_i(\sigma)}\otimes A_{io}\right)\varepsilon_{i\star}\notag\\
		&-\left(\mathcal{H}_{\pi_i(\sigma)}\otimes I_{v_i}\right)\left(\Gamma_i^{\sigma}\otimes I_{v_i}\right)\left(\mathcal{H}_{\pi_i(\sigma)}\otimes I_{v_i}\right)\varepsilon_{i\star}\notag\\
		&+\left(\mathcal{H}_{\pi_i(\sigma)}\otimes I_{v_i}\right)\left(\Gamma_i^{\sigma}\otimes I_{v_i}\right)\left(B_{i\star}^{\sigma}\otimes I_{v_i}\right)\varepsilon_{ii}\notag\\
		&-\left(B_{i\star}^{\sigma}\otimes I_{v_i}\right)\left(I_{\pi_i(\sigma)}\otimes \left(A_{io}-H_{io}C_{io}\right)\right)\varepsilon_{ii}\notag\\
		&-\sum_{l=1}^{i-1}\left(B_{i\star}^{\sigma}\otimes I_{v_i}\right)\left(I_{\pi_i(\sigma)}\otimes \Upsilon_{il}\right)\left(1_{\pi_i(\sigma)}\otimes e_{i,l}\right)\notag\\
		&+\sum_{l=1}^{i-1}\left(\mathcal{H}_{\pi_i(\sigma)}\otimes I_{v_i}\right)\Upsilon_{i\star,l}\varepsilon_{i\star,l}.	
	\end{align}
	According to Lemma \ref{renet}, subgraph $\mathcal{V}_i^\sigma$ has a spanning tree rooted at $i$. Therefore, $\mathcal{H}_{\pi_i(\sigma)}$ is full rank due to the conclusion of Lemma \ref{spanning}. Then, we have
	\begin{align}\label{istar}
		\varepsilon_{i\star}=&\left(\mathcal{H}_{\pi_i(\sigma)}\otimes I_{v_i}\right)^{-1}\xi_{i\star}\notag\\
		&+\left(\mathcal{H}_{\pi_i(\sigma)}\otimes I_{v_i}\right)^{-1}\left(B_{i\star}^{\sigma}\otimes I_{v_i}\right)\varepsilon_{ii}.
	\end{align}
	Subsequently, substituting (\ref{istar}) into (\ref{xi}) leads to
	\begin{align*}
		\dot{\xi}_{i\star}=&\left(I_{\pi_i(\sigma)}\otimes A_{io}\right)\xi_{i\star}+\left(B_{i\star}^{\sigma}\otimes A_{io}\right)\varepsilon_{ii}\notag\\
		&-\left(\mathcal{H}_{\pi_i(\sigma)}\otimes I_{v_i}\right)\left(\Gamma_i^{\sigma}\otimes I_{v_i}\right)\xi_{i\star}\notag\\
		&-\left(B_{i\star}^{\sigma}\otimes I_{v_i}\right)\left(I_{\pi_i(\sigma)}\otimes \left(A_{io}-H_{io}C_{io}\right)\right)\varepsilon_{ii}\notag\\
		&-\sum_{l=1}^{i-1}\left(B_{i\star}^{\sigma}\otimes I_{v_i}\right)\left(I_{\pi_i(\sigma)}\otimes \Upsilon_{il}\right)\left(1_{\pi_i(\sigma)}\otimes e_{i,l}\right)\notag\\
		&+\sum_{l=1}^{i-1}\left(\mathcal{H}_{\pi_i(\sigma)}\otimes I_{v_i}\right)\Upsilon_{i\star,l}\varepsilon_{i\star,l}\notag\\
		=&\left(I_{\pi_i(\sigma)}\otimes A_{io}\right)\xi_{i\star}-\left(\mathcal{H}_{\pi_i(\sigma)}\Gamma_i^{\sigma}\otimes I_{v_i}\right)\xi_{i\star}\notag\\
		&+\left(B_{i\star}^{\sigma}\otimes H_{io}C_{io}\right)\varepsilon_{ii}\notag\\
		&-\sum_{l=1}^{i-1}\left(B_{i\star}^{\sigma}\otimes I_{v_i}\right)\left(I_{\pi_i(\sigma)}\otimes \Upsilon_{il}\right)\left(1_{\pi_i(\sigma)}\otimes e_{i,l}\right)\notag\\
		&+\sum_{l=1}^{i-1}\left(\mathcal{H}_{\pi_i(\sigma)}\otimes I_{v_i}\right)\Upsilon_{i\star,l}\varepsilon_{i\star,l}.
	\end{align*}
	Hence, (\ref{xi2}) is obtained.
\end{IEEEproof}

\begin{rem}
	This remark describes the significance of observability decomposition and state reorganization methods. The error dynamics of each local observer contains the dynamics of $N$ observable subsystems, and the convergence of the error dynamics corresponding to different observable subsystems depends on different network topology. Therefore, the communication conditions can meet the requirements of all observable subsystems in all local observers at the same time if and only if the network is strongly connected. However, in the scenario of jointly connected switching topologies, the real-time disconnected network and the strong coupling relationship between observable subsystems will make it impossible to analyze the performance of the error dynamics of any observable subsystem. Hence, the observability decomposition method is employed to decompose the state estimation of each local observer, which leads to a weak coupling relationship between all decomposed systems. Furthermore, according to the communication network requirements of the decomposed observable subsystems, the subsystems from each local observer with the same communication requirements are reorganized into a group of new compact states including $\varepsilon_{ii}, \varepsilon_{i\star}$, and $\varepsilon_{i\diamond}$. In this way, all state components in new compact states can meet the communication conditions required for stability analysis. Consequently, the network properties of each new compact state and the weak coupling among them help us construct the recursive proof method in the next section, so as to obtain the asymptotic omniscience of the distributed observer under the switching topologies.
\end{rem}

\begin{rem}\label{Xi}
	As stated in the previous remark, each new compact state has a complete network property that can meet its own requirements after state decomposition and reorganization. Herein, we elaborate on the specific meaning of the complete network property. Denote $\Xi_{ij}$ the $j$th observable subsystem of the $i$th local observer after observability decomposition. A subsystem $\Xi_{ij}$ has a complete network property if all of its parents and children nodes are descendants of $\Xi_{jj}$ when $\Xi_{ij}$ is a descendant node of $\Xi_{jj}$ (note that $\Xi_{ij}$ is a descendant node of $\Xi_{jj}$ if and only if $i$ is the descendant of $j$ because $i$ in $\Xi_{ij}$ labels the local observer). It is impossible for a typical network topology to possess such properties, and the key to ensuring that the communication network maintains such properties during the switching process in this paper lies in the network transformation mapping proposed in Section \ref{sec3}. The network transformation mapping with respect to node $j$ can ensure that the descendant nodes of $j$ will not receive information from other non-descendant nodes. Therefore, it can ensure that a node $i$ is either not a descendant of $j$ or it and its parents and children are all descendants of $j$ (a more detailed description can be found in Remark \ref{rem-ntm}). In summary, the state decomposition and reorganization method in this section and the network transformation mapping in the previous section are important foundations for the asymptotic omniscience proof in the next section.
\end{rem}

\subsection{Stability analysis of the error dynamics of distributed observer}\label{sec4.2}
By observing the form (\ref{xi2}) of $\xi_{i\star}$, it is not difficult to find that there is a lower triangular coupling relationship between $\xi_{1\star},\ldots,\xi_{N\star}$. Namely, the dynamics of $\xi_{i\star}$ are affected by $e_{i,t}$ with $t=1,\ldots,i-1$. Their stability is proved recursively in light of this coupling relationship, which is demonstrated in the following theorem.

\begin{thm}\label{thm1}
	Distributed observer formed by (\ref{do1}), (\ref{do2}) achieves omniscience asymptotically under the switching directed networks subject to Assumptions \ref{assume3} and \ref{assume2}, if its coupling gain satisfies (\ref{doubler}) with a positive constant $\gamma$ and a group of sufficiently large $\gamma_{i,k}$ governed by 
	\begin{align}
		\underline{\gamma}_i^2&>\bar{\lambda}(\Xi_{i1})+\frac{1}{4}\bar{\lambda}^2(\Xi_{i2})+\wp,\label{gammai}\\
		\wp&>\frac{\bar{\lambda}_i\bar{\lambda}(\Xi_{i\diamond})}{\underline{\lambda}(P_{iu})}\frac{m(T_{k,i}^c)}{m(T_{k,i})},\label{W_i}
	\end{align}
	for all $T_{k,i}$ and $T_{k,i}^c$ with $i=1,\ldots,N$ and $k=1,\ldots,N$. Here, $m(\cdot)$ is measure and $T_{k,i}$ is a subset of $[t_j,t_{j+1})$ in which $k\in\mathcal{V}_{i}^{\sigma}$, and $T_{k,i}^c$ is $[t_j,t_{j+1})\backslash T_{k,i}$; 
	$\underline{\gamma}_i=\min_{j=1,\ldots,N}\{\gamma_{j,i}\}$;
	$\Xi_{i1}^\sigma=Q_{io}^{\sigma}\otimes A_{io}$, $\Xi_{i2}^\sigma=Q_{io}^{\sigma}B_{i\star}^{\sigma}\otimes H_{io}C_{io}$, $\bar{\lambda}(\Xi_{i1})=\max_{t\in\mathbb{R}}\{\bar{\lambda}(\Xi_{i1}^\sigma)\}$, $\bar{\lambda}(\Xi_{i2})=\max_{t\in\mathbb{R}}\{\bar{\lambda}(\Xi_{i2}^\sigma)\}$; 
	$\wp>0$ is a given constant; 
	$Q_{io}^\sigma$ satisfies $Q_{io}^{\sigma}\mathcal{H}_{\pi_i(\sigma)}+\mathcal{H}_{\pi_i(\sigma)}^TQ_{io}^{\sigma}>2I_{\pi_i(\sigma)}$; $\bar{\lambda}_i^\sigma=\max\left\{\bar{\lambda}(P_{io}),\bar{\lambda}\left(Q_{io}^{\sigma}\right)\right\}$ and $\bar{\lambda}_i=\max_{t\in\mathbb{R}}\{\bar{\lambda}_i^\sigma\}$; $\Xi_{i\diamond}^\sigma=2P_{iu}^\sigma\otimes A_{1o}$ with symmetric positive definite matrix $P_{iu}^\sigma$ being solved by $P_{iu}^\sigma\mathcal{L}_{\bar{\pi}_i(\sigma)}+\mathcal{L}_{\bar{\pi}_i(\sigma)}^TP_{iu}^\sigma\geq 0$, and $\mathcal{L}_{\bar{\pi}_1(\sigma)}$ is the Laplacian matrix that corresponds to the subgraph among $\mathcal{V}\backslash\mathcal{V}_i^\sigma$.
\end{thm}

\begin{IEEEproof}
	This theorem will be proved by mathematical induction. 
	
	1)~The performance of $\varepsilon_{11}$ and $\xi_{1\star}$. We will show that $\lim_{t\to\infty}\|\varepsilon_{11}\|=0$ and $\lim_{t\to\infty}\|\xi_{1\star}\|=0$. According to (\ref{varepsilonii}), it has
	\begin{align*}
		\dot{\varepsilon}_{11}=\left(I_{\pi_1(\sigma)}\otimes (A_{1o}-H_{1o}C_{1o})\right)\varepsilon_{11}.
	\end{align*}
	A symmetric positive definite matrix $P_{1o}$ can be found such that $sym\left\{P_{1o}(A_{1o}-H_{1o}C_{1o})\right\}=-2(\bar{\gamma}_1+\wp)I_{v_1}$,
	where $\bar{\gamma}_1=\max_{j=1,\ldots,N}\{\gamma_{j,1}\}$. Then, by setting $V_{1,\varepsilon}=\frac{1}{2}\varepsilon_{11}^T(I_{\pi_1(\sigma)}\otimes P_{1o})\varepsilon_{11}$, it yields $\dot{V}_{1,\varepsilon}\leq -(\bar{\gamma}_1+\wp)\varepsilon_{11}^T\varepsilon_{11}$. Also, since $\mathcal{H}_{\pi_1(\sigma)}$ is a nonsingular M-matrix, there exists positive diagonal matrix $Q_{1o}^{\sigma}$ that satisfies $Q_{1o}^{\sigma}\mathcal{H}_{\pi_1(\sigma)}+\mathcal{H}_{\pi_1(\sigma)}^TQ_{1o}^{\sigma}>2I_{\pi_1(\sigma)}$.
	To move on, the Lyapunov function with respect to $\xi_{1\star}$ can be chosen as 
	\begin{align*}
		V_{1,\xi}=\frac{1}{2}\sum_{j\in\mathcal{V}_{1}^{\sigma}\backslash\{1\}}q_{1o,j}^{\sigma}\gamma_{j,1}\xi_{j,1}^T\xi_{j,1},
	\end{align*}
	with $q_{1o,j}^\sigma$ being the elements of a symmetric positive definite matrix $Q_{1o}^{\sigma}=diag\{q_{1o,j},~j\in\mathcal{V}_{1}^{\sigma}\backslash\{1\}\}$.
	It follows that
	\begin{align*}
		\dot{V}_{1,\xi}=&\xi_{1\star}^T\left(\Gamma_1^{\sigma}Q_{1o}^{\sigma}\otimes I_{v_1}\right)\dot{\xi}_{1\star}\notag\\
		=&\xi_{1\star}^T\left(\Gamma_1^{\sigma}Q_{1o}^{\sigma}\otimes A_{1o}-\Gamma_1^{\sigma}Q_{1o}^{\sigma}\mathcal{H}_{\pi_1(\sigma)}\Gamma_1^{\sigma}\otimes I_{v_1}\right)\xi_{1\star}\notag\\
		&+\xi_{1\star}^T\left(\Gamma_1^{\sigma}Q_{1o}^{\sigma}B_{1\star}^{\sigma}\otimes H_{1o}C_{1o}\right)\varepsilon_{11}\notag\\
		\leq& \bar{\lambda}\left(Q_{1o}^{\sigma}\otimes A_{1o}\right)\xi_{1\star}^T\left(\Gamma_1^{\sigma}\otimes I_{v_1}\right)\xi_{1\star}\notag\\
		&-\xi_{1\star}^T\left(\Gamma_1^{\sigma}\Gamma_1^{\sigma}\otimes I_{v_1}\right)\xi_{1\star}\notag\\
		&+\bar{\lambda}\left(Q_{1o}^{\sigma}B_{1\star}^{\sigma}\otimes H_{1o}C_{1o}\right)\xi_{1\star}^T\left(\Gamma_1^{\sigma}\otimes I_{v_1}\right)\varepsilon_{11}\notag\\
		\leq& \left(\bar{\lambda}(Q_{1o}^{\sigma}\otimes A_{1o})-\underline{\gamma}_{1}^2\right)\xi_{1\star}^T\left(\Gamma_1^{\sigma}\otimes I_{v_1}\right)\xi_{1\star}\notag\\
		&+\bar{\lambda}\left(Q_{1o}^{\sigma}B_{1\star}^{\sigma}\otimes H_{1o}C_{1o}\right)\xi_{1\star}^T\left(\Gamma_1^{\sigma}\otimes I_{v_1}\right)\varepsilon_{11},
	\end{align*}
	where $\underline{\gamma}_1=\min_{j=1,\ldots,N}\{\gamma_{j,1}\}$.
	
	For ease of notation, we define $\wp_{1,1}=-\bar{\gamma}_1+\wp$, $\wp_{1,2}=\bar{\lambda}\left(Q_{1o}^{\sigma}B_{1\star}^{\sigma}\otimes H_{1o}C_{1o}\right)$, $\wp_{1,3}=\bar{\lambda}(Q_{1o}^{\sigma}\otimes A_{1o})-\underline{\gamma}_{1}^2+\wp$, and $\bar{\Gamma}_1^{\sigma}=\Gamma_1^{\sigma}\otimes I_{v_1}$. Let $V_{1,1}=V_{1,\varepsilon}+V_{1,\xi}$ and hence deduce
	\begin{align}\label{V1}
		\dot{V}_{1,1}\leq& \wp_{1,3}\xi_{1\star}^T\bar{\Gamma}_1^{\sigma}\xi_{1\star}+\wp_{1,2}\xi_{1\star}^T\bar{\Gamma}_1^{\sigma}\varepsilon_{11}-\bar{\gamma}_1\varepsilon_{11}^T\varepsilon_{11}\notag\\
		&-\wp\xi_{1\star}^T\bar{\Gamma}_1^{\sigma}\xi_{1\star}-\wp\varepsilon_{11}^T\varepsilon_{11}\notag\\
		=&\begin{bmatrix}\varepsilon_{11}^T&\xi_{1\star}^T\end{bmatrix}\begin{bmatrix}-\bar{\gamma}_1I_{\pi_{1}(\sigma)v_1}&\frac{1}{2}\wp_{1,2}\bar{\Gamma}_1^{\sigma}\\ \frac{1}{2}\wp_{1,2}\bar{\Gamma}_1^{\sigma}&\wp_{1,3}\bar{\Gamma}_1^{\sigma}\end{bmatrix}\begin{bmatrix}\varepsilon_{11}\\ \xi_{1\star}\end{bmatrix}\notag\\
		&-\begin{bmatrix}\varepsilon_{11}^T&\xi_{1\star}^T\end{bmatrix}\begin{bmatrix}\wp I_{\pi_{1}(\sigma)v_1}&\\ &\wp\bar{\Gamma}_1^{\sigma}\end{bmatrix}\begin{bmatrix}\varepsilon_{11}\\ \xi_{1\star}\end{bmatrix}.
	\end{align}
	Based on (\ref{gammai}), we know
	\begin{align}\label{gamma_bar}
		\underline{\gamma}_1^2&>\bar{\lambda}(\Xi_{11})+\frac{1}{4}\bar{\lambda}^2(\Xi_{12})+\wp\notag\\
		&>\bar{\lambda}(\Xi_{11}^\sigma)+\frac{1}{4}\bar{\lambda}^2(\Xi_{12}^\sigma)+\wp,
	\end{align}
	where $\Xi_{11}^\sigma=Q_{1o}^{\sigma}\otimes A_{1o}$, $\Xi_{12}^\sigma=Q_{1o}^{\sigma}B_{1\star}^{\sigma}\otimes H_{1o}C_{1o}$, $\bar{\lambda}(\Xi_{11})=\max_{t\in\mathbb{R}}\{\bar{\lambda}(\Xi_{11}^\sigma)\}$, $\bar{\lambda}(\Xi_{12})=\max_{t\in\mathbb{R}}\{\bar{\lambda}(\Xi_{12}^\sigma)\}$. Hence,
	\begin{align*}
		0&>\wp_{1,3}I_{\pi_1(\sigma)v_1}+\frac{1}{4}\wp_{1,2}^2I_{\pi_1(\sigma)v_1}\notag\\
		&>\wp_{1,3}I_{\pi_1(\sigma)v_1}+\frac{1}{4\bar{\gamma}_1}\wp_{1,2}^2\left(\Gamma_1^{\sigma}\otimes I_{v_1}\right).
	\end{align*}
	Since $\Gamma_1^{\sigma}\otimes I_{v_1}$ is a diagonal positive definite matrix, we further obtain that
	\begin{align}\label{schur1}
		\wp_{1,3}\left(\Gamma_1^{\sigma}\otimes I_{v_1}\right)+\frac{1}{4\bar{\gamma}_1}\wp_{1,2}^2\left(\Gamma_1^{\sigma}\otimes I_{v_1}\right)^2<0.
	\end{align}
	In light of the Schur complement lemma, (\ref{schur1}) is equivalent to
	\begin{align}\label{schur2}
		\begin{bmatrix}-\bar{\gamma}_1I_{\pi_{1}(\sigma)v_1}&\frac{1}{2}\wp_{1,2}\bar{\Gamma}_1^{\sigma}\\ \frac{1}{2}\wp_{1,2}\bar{\Gamma}_1^{\sigma}&\wp_{1,3}\bar{\Gamma}_1^{\sigma}\end{bmatrix}<0.
	\end{align}
	
	Substituting (\ref{schur2}) into (\ref{V1}) yields
	\begin{align*}
		\dot{V}_{1,1}&<-\wp\varepsilon_{11}^T\varepsilon_{11}-\wp\xi_{1\star}^T\bar{\Gamma}_1^{\sigma}\xi_{1\star}\notag\\
		&<-\frac{\wp}{\bar{\lambda}(P_{1o})}V_{1,\varepsilon}-\frac{\wp}{\bar{\lambda}\left(Q_{1o}^{\sigma}\right)}V_{1,\xi}\notag\\
		&<-(\bar{\lambda}_1^\sigma)^{-1}\wp V_{1,1},
	\end{align*}
	where $\bar{\lambda}_1^\sigma=\max\left\{\bar{\lambda}(P_{1o}),\bar{\lambda}\left(Q_{1o}^{\sigma}\right)\right\}$. Since the range of $\sigma(t)$ is limited, there exists $\bar{\lambda}_1=\max_{t\in\mathbb{R}}\{\bar{\lambda}_1^\sigma\}$. Hence 
	\begin{align*}
		V_{1,1}(t)\leq \exp\left\{-\bar{\lambda}_1^{-1}\wp(t-t_0)\right\}V_{1,1}(t_0).
	\end{align*}
	
	Now the nodes that belong to $\mathcal{V}\backslash\mathcal{V}_{1}^{\sigma}$ should be considered. Note that $\ell_{ik}^{\sigma}=0$ in this scenario and thus $\mathbbm{r}_{i,k}=\gamma$. Recall that the compact form of these nodes is $\varepsilon_{1\diamond}=col\{e_{k,1},~k\in\mathcal{V}\backslash\mathcal{V}_{1}^{\sigma}\}$. By setting $\bar{\pi}_{1}(\sigma)=N-\pi_1(\sigma)$, we obtain
	\begin{align*}
		\dot{\varepsilon}_{1\diamond}=\left(I_{\bar{\pi}_{1}(\sigma)}\otimes A_{1o}\right)\varepsilon_{1\diamond}-\gamma\left(\mathcal{L}_{\bar{\pi}_{1}(\sigma)}\otimes I_{v_1}\right)\varepsilon_{1\diamond}.
	\end{align*}
	Since $\mathcal{L}_{\bar{\pi}_1(\sigma)}$ is a singular M-matrix, there exists a symmetric positive definite matrix $P_{1u}^\sigma$ in light of \cite[Chapter 4]{Cooperative2009Qu} such that $P_{1u}^\sigma\mathcal{L}_{\bar{\pi}_1(\sigma)}+\mathcal{L}_{\bar{\pi}_1(\sigma)}^TP_{1u}^\sigma$ is a symmetric positive semidefinite matrix. Subsequently, the Lyapunov candidate can be chosen as $V_{1,\diamond}(t)=\varepsilon_{1\diamond}^T\left(P_{1u}^\sigma\otimes I_{\bar{\pi}_{1}(\sigma)}\right)\varepsilon_{1\diamond}$, which leads to
	\begin{align*}
		\dot{V}_{1,\diamond}(t)=&\varepsilon_{1\diamond}^T\left(P_{1u}^\sigma\otimes A_{1o}\right)\varepsilon_{1\diamond}+\varepsilon_{1\diamond}^T\left(P_{1u}^\sigma\otimes A_{1o}^T\right)\varepsilon_{1\diamond}\notag\\
		&-\gamma\varepsilon_{1\diamond}^T\left(\left(P_{1u}^\sigma\mathcal{L}_{\bar{\pi}_1(\sigma)}+\mathcal{L}_{\bar{\pi}_1(\sigma)}^TP_{1u}^\sigma\right)\otimes I_{v_1}\right)\varepsilon_{1\diamond}\notag\\
		\leq& 2\varepsilon_{1\diamond}^T\left(P_{1u}^\sigma\otimes A_{1o}\right)\varepsilon_{1\diamond}\leq\frac{\bar{\lambda}(\Xi_{1\diamond}^\sigma)}{\underline{\lambda}(P_{1u}^\sigma)}V_{1,\diamond}(t),
	\end{align*}
	where $\Xi_{1\diamond}^\sigma=2P_{1u}^\sigma\otimes A_{1o}$. Also, because the range of $\sigma$ is finite, we find $\bar{\lambda}(\Xi_{1\diamond})=\max_{t\in\mathbb{R}}\{\bar{\lambda}(\Xi_{1\diamond}^\sigma)\}$, $\underline{\lambda}(P_{1u})=\min_{t\in\mathbb{R}}\{\underline{\lambda}(P_{1u}^\sigma)\}$, and
	\begin{align*}
		V_{1,\diamond}(t)\leq\exp\left\{\frac{\bar{\lambda}(\Xi_{1\diamond})}{\underline{\lambda}(P_{1u})}(t-t_0)\right\}V_{1,\diamond}(t_0).
	\end{align*}
	
	Let $W_{i,1}$ be the Lyapunov function of $e_{i,1}$. It satisfies
	\begin{align*}
		W_{i,1}(t)=
		e_{1,1}^TP_{1o}e_{1,1}+\frac{1}{2}\underline{\lambda}(Q_{1o})\gamma_{i,1}\xi_{i,1}^T\xi_{i,1}
	\end{align*}
	if $t\in T_{i,1}$, and
	\begin{align*}
		W_{i,1}(t)=\underline{\lambda}(P_{1u})e_{i,1}^Te_{i,1}
	\end{align*}
	if $t\in T_{i,1}^c$, where $T_{i,1}$ is a subset of $[t_k,t_{k+1})$ in which $i\in\mathcal{V}_{1}^{\sigma}$, and $T_{i,1}^c$ is $[t_k,t_{k+1})\backslash T_{i,1}$. Let $\beth_{i,1}^{k\chi}$ be the convergence (divergence) rate of $W_{i,1}(t)$ when $t\in [t_k^{\chi-1},t_k^\chi)$. Noting that $\sum_{i\in\mathcal{V}_1^\sigma\backslash\{1\}}W_{i,1}(t)\leq V_{1,1}$ if $[t_k^{\chi-1},t_k^\chi)\subset T_{i,1}$, and $\sum_{i\in\mathcal{V}\backslash\mathcal{V}_1^\sigma}W_{i,1}(t)\leq V_{1,\diamond}$ if $[t_k^{\chi-1},t_k^\chi)\subset T_{i,1}^c$. Therefore, $\beth_{i,1}^{k\chi}\leq-\bar{\lambda}_1^{-1}\wp$ if $[t_k^{\chi-1},t_k^\chi)\subset T_{i,1}$ and $\beth_{i,1}^{k\chi}\leq\frac{\bar{\lambda}(\Xi_{1\diamond})}{\underline{\lambda}(P_{1u})}$ otherwise. Then,
	\begin{align*}
		&W_{i,1}(t_{k+1})=W_{i,1}\left(t_k^\ell\right)\notag\\
		=& \exp\left\{\beth_{i,1}^{k\ell}(t_k-t_k^{\ell-1})\right\}W_{i,1}\left(t_k^{\ell-1}\right)\notag\\
		\leq& \exp\left\{\beth_{i,1}^{k\ell}(t_k^\ell-t_k^{\ell-1})\right\}\exp\left\{\beth_{i,1}^{k(\ell-1)}(t_k^{\ell-1}-t_k^{\ell-2})\right\}\notag\\
		&\times W_{i,1}\left(t_k^{\ell-2}\right)\notag\\
		\leq& \prod_{\chi=1}^\ell\exp\left\{\beth_{i,1}^{k\chi}(t_k^\chi-t_k^{\chi-1})\right\}W_{i,1}\left(t_k^{0}\right).
	\end{align*}
	According to the value of $\beth_{i,1}^{k\chi}$, we know
	\begin{align*}
		&\prod_{\chi=1}^\ell\exp\left\{\beth_{i,1}^{k\chi}(t_k^\chi-t_k^{\chi-1})\right\}\notag\\
		=&\exp\left\{\sum_{\chi=1}^\ell\beth_{i,1}^{k\chi}(t_k^\chi-t_k^{\chi-1})\right\}\notag\\
		=&\exp\left\{-\bar{\lambda}_1^{-1}\wp m(T_{i,1})+\frac{\bar{\lambda}(\Xi_{1\diamond})}{\underline{\lambda}(P_{1u})}m(T_{i,1}^c)\right\}.
	\end{align*}
	Consequently, 
	\begin{align*}
		W_{i,1}(t_{k+1})\leq\exp\left\{\beth_{i,1}\right\}W_{i,1}(t_k)
	\end{align*}
	with
	\begin{align*}
		\beth_{i,1}=-\bar{\lambda}_1^{-1}\wp m(T_{i,1})+\frac{\bar{\lambda}(\Xi_{1\diamond})}{\underline{\lambda}(P_{1u})}m(T_{i,1}^c).
	\end{align*}
	(\ref{W_i}) indicates that $W_{i,1}(t_{k+1})\leq W_{i,1}(t_k)$ because it leads to $\beth_{i,1}<0$. Therefore, the sufficient conditions (\ref{gammai}) and (\ref{W_i}) guarantee the convergence of $\varepsilon_{11}$ and $\xi_{1\star}$.
	
	2)~Suppose the convergence of $\varepsilon_{kk}$ and $\xi_{k\star}$ are fulfilled for $k=1,\ldots,i-1$. We will show that $\lim_{t\to\infty}\|\varepsilon_{ii}\|=0$ and $\lim_{t\to\infty}\|\xi_{i\star}\|=0$. From Lemma \ref{lem7}, we note 
	\begin{align}\label{varepsilonii2}
		\dot{\varepsilon}_{ii}=&\left(I_{\pi_i(\sigma)}\otimes \left(A_{io}-H_{io}C_{io}\right)\right)\varepsilon_{ii}\notag\\
		&+\sum_{l=1}^{i-1}\left(I_{\pi_i(\sigma)}\otimes \Upsilon_{il}\right)\left(1_{\pi_i(\sigma)}\otimes e_{i,l}\right)
	\end{align}
	and
	\begin{align}
		\dot{\xi}_{i\star}=&\left(I_{\pi_i(\sigma)}\otimes A_{io}\right)\xi_{i\star}-\left(\mathcal{H}_{\pi_i(\sigma)}\Gamma_i^{\sigma}\otimes I_{v_i}\right)\xi_{i\star}\notag\\
		&+\left(B_{i\star}^{\sigma}\otimes H_{io}C_{io}\right)\varepsilon_{ii}\notag\\
		&-\sum_{l=1}^{i-1}\left(B_{i\star}^{\sigma}\otimes I_{v_i}\right)\left(I_{\pi_i(\sigma)}\otimes \Upsilon_{it}\right)\left(1_{\pi_i(\sigma)}\otimes e_{i,l}\right)\notag\\
		&+\sum_{l=1}^{i-1}\left(\mathcal{H}_{\pi_i(\sigma)}\otimes I_{v_i}\right)\Upsilon_{i\star,l}\varepsilon_{i\star,l}.
	\end{align}
	Since $\lim_{t\to\infty}\|e_{i,t}\|=0$, $\lim_{t\to\infty}\|\varepsilon_{i\star,t}\|=0$ for $i=1,\ldots,t-1$, in light of Lemma \ref{cyb}, the stability of $\varepsilon_{ii}$ and $\xi_{i\star}$ is equivalent to the stability
	\begin{align}
		&\dot{\epsilon}_{ii}=\left(I_{\pi_i(\sigma)}\otimes (A_{io}-H_{io}C_{io})\right)\epsilon_{ii},\label{sysiio}\\
		&\dot{\zeta}_{i\star}=\left(I_{\pi_i(\sigma)}\otimes A_{io}\right)\zeta_{i\star}-\left(\mathcal{H}_{\pi_i(\sigma)}\Gamma_i^{\sigma}\otimes I_{v_i}\right)\zeta_{i\star}\notag\\
		&~~~~~~~~+\left(B_{i\star}^{\sigma}\otimes H_{io}C_{io}\right)\epsilon_{ii},\label{sysipo}
	\end{align}   
	where the newly defined notation $\epsilon_{ii}$ and $\zeta_{i\star}$ are employed to replace $\varepsilon_{ii}$ and $\xi_{i\star}$ to avoid confusion. By repeating the same process in the stability proof of the first subsystem, we can obtain that the states (\ref{sysiio}), (\ref{sysipo}) converge to zero if (\ref{gammai}) and (\ref{W_i}) hold.
	
	Therefore, the trajectories of (\ref{varepsilonii}) and (\ref{xi2}) converge to zero exponentially; so do $\varepsilon_{ii}$ and $\varepsilon_{i\star}$, owing to Lemma \ref{lem7}. This implies the convergence of all $e_i$ with $i=1,\ldots,N$.
\end{IEEEproof}
\begin{rem}
	\color{red} We can see from the proof of Theorem \ref{thm1} that the stability of the error dynamics of the distributed observer under the switching topology designed in this paper does not depend on the initial value, which is also consistent with general understanding of linear system theory.
\end{rem}
\begin{rem}
	This section proves a sufficient condition for a distributed observer to achieve asymptotic omniscience in the case of directed switching networks. Distributed observer dynamics are decomposed using the network transformation mapping proposed in Section \ref{sec3}, and the decomposed observer dynamics can make full use of the properties of the directed support tree. 
\end{rem}

\section{\color{pinegreen}Achievement of distributed observer with adaptive coupling gain}\label{sec5}

{\color{pinegreen}Due to the complex calculation of the sufficient conditions (\ref{gammai}), adaptive coupling gain will be proposed in this section.}

One of the main roles of the unique design of $\mathbbm{r}_{i,k}$ in (\ref{do1}), (\ref{do2}) is that it helps us to convert the constant coupling gain into an adaptive one. In this paper, the adaptive gain is only employed in the subsystems $\Xi_{ik}$ (see the definition in Remark \ref{Xi}) that satisfy $\ell_{i,k}^\sigma=1$, because it is unnecessary to adaptively adjust the coupling gain in subsystems without a convergence trend. For this reason, the adaptive gain is introduced as
\begin{align}
	&\dot{\gamma}_{i,k}=\omega_{i,k},~\gamma_{i,k}(0)>0,\label{ada-gamma}\\
	&\omega_{i,k}=\left\|\sum_{j=1}^N\alpha_{ij}^{\sigma(t)}(k)\left(\hat{x}_{j,k}-\hat{x}_{i,k}\right)\right\|^2.
\end{align}

The corresponding Lyapunov function chosen is $V_{i,1}^a=V_{i,\varepsilon}^a+V_{i,\xi}^a$, where
\begin{align}
	V^a_{i,\varepsilon}=&\frac{1}{2}\varepsilon_{ii}^T(I_{\pi_i(\sigma)}\otimes P_{io})\varepsilon_{ii},\\
	V^a_{i,\xi}=&\frac{1}{2}\sum_{j\in\mathcal{V}_{i}^{\sigma(t)}\backslash\{i\}}\left(q_{io,j}^\sigma(2\gamma_{j,i}+\omega_{j,i})\omega_{j,i}\right.\notag\\
	&\left.+q_{io,j}^\sigma(\gamma_{j,i}-\gamma_i^*)^2\right),
\end{align}
with a pending decision constant $\gamma_i^*$ and positive constant $q_{io,j}^\sigma$ that belongs to $Q_{io}^\sigma=col\{q_{io,j}^\sigma,~j\in\mathcal{V}_{i}^{\sigma}\backslash\{i\}\}$. Then, we can directly show that $\varepsilon_{ii}$ and $\xi_{i\star}$ have a convergence trend (see the corresponding proof in the appendix). However, with the adaptive gain, we cannot estimate the convergence rate when the node satisfies $\ell_{i,k}^\sigma=1$. Hence, we fail to guarantee that the error dynamics of $\hat{x}_{i,k}$ converge to zero with $t\to\infty$.

Fortunately, the proof in Theorem \ref{thm1} shows that the convergence rate of $V_{k,k}$ depends on the coupling gain $\gamma_{i,k},~i\in\mathcal{V}_{k}^\sigma$. This ensures that the error dynamics can converge at least at the speed of $\wp$ when (\ref{gammai}) is satisfied. Actually, the overall error dynamics cannot converge if the error convergence rate of subsystems on those nodes satisfying $\ell_{i,k}^\sigma$ is less than $\wp$. It thus leads to the increasing value of $\gamma_{i,k}$, until the coupling gain fulfills the requirement of (\ref{gammai}). The above statement indicates that the distributed observer can achieve omniscience asymptotically with adaptive gain under the switching directed networks.
\begin{figure}[!t]
	\centering
	\includegraphics[width=6.5cm]{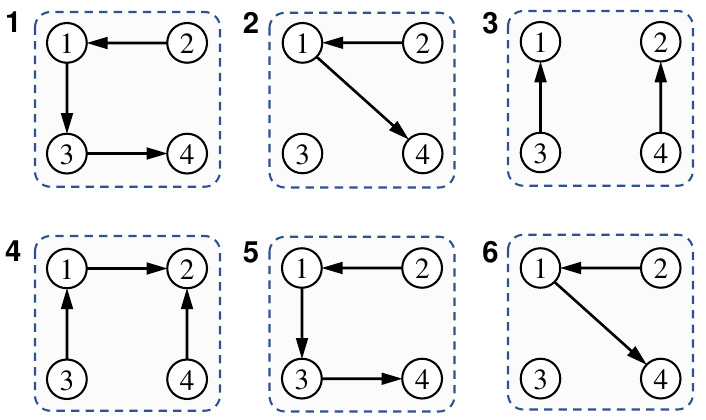}\\
	\caption{Augmented communication networks with $4$ agents}\label{topology}
\end{figure}

\begin{table}\caption{Parameters of the power system}\label{table}
	\renewcommand\arraystretch{1.2}
	\begin{tabular}{p{1cm}<{\centering}p{1.3cm}<{\centering}p{1.3cm}<{\centering}p{1.3cm}<{\centering}p{1.3cm}<{\centering}}
		\hline
		~     &Area 1   & Area 2  & Area 3 & Area 4 \\
		\hline
		$M_i$    & $12$    & $10$    & $8$    & $8$ \\
		$R_i$    & $0.05$  & $0.0625$& $0.08$ & $0.08$\\
		$D_i$    & $0.7$   & $0.9$   & $0.9$  & $0.7$ \\
		$T_{t_i}$& $0.65$  & $0.4$   & $0.3$  & $0.6$ \\
		$T_{g_i}$& $0.1$   & $0.1$   & $0.1$  & $0.1$ \\
		\hline    
	\end{tabular}
\end{table}

\begin{figure*}[!t]
	\centering
	\includegraphics[width=18cm]{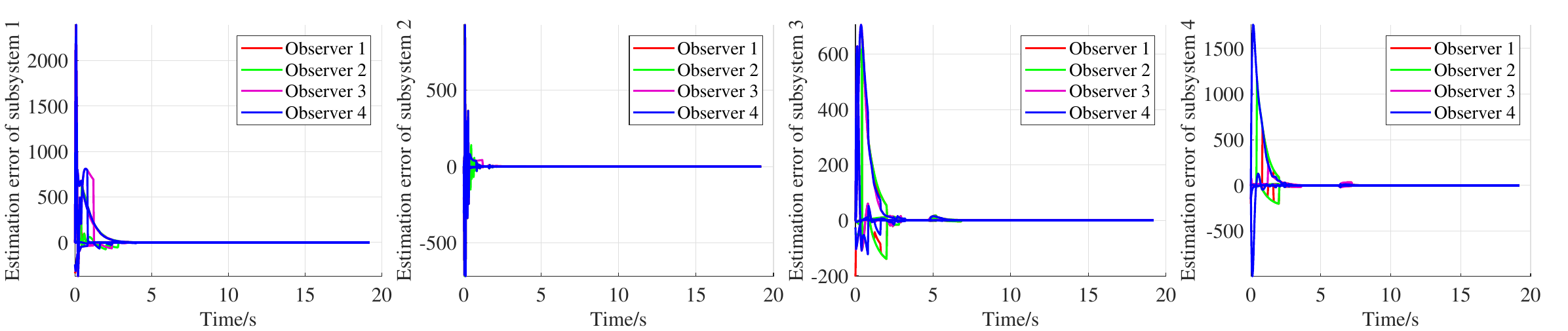}\\
	\caption{Error dynamics of distributed observer with network transformation mapping}\label{fig-powererror}
\end{figure*}

\begin{figure*}[!t]
	\centering
	\includegraphics[width=18cm]{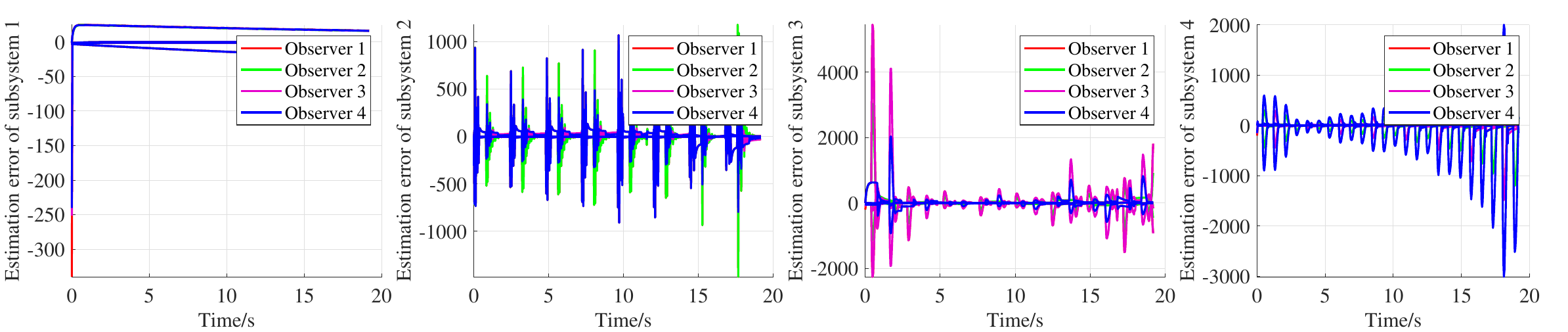}\\
	\caption{Error dynamics of distributed observer without network transformation mapping}\label{fig-powerwithout}
\end{figure*}

\begin{figure*}[!t]
	\centering
	\includegraphics[width=18cm]{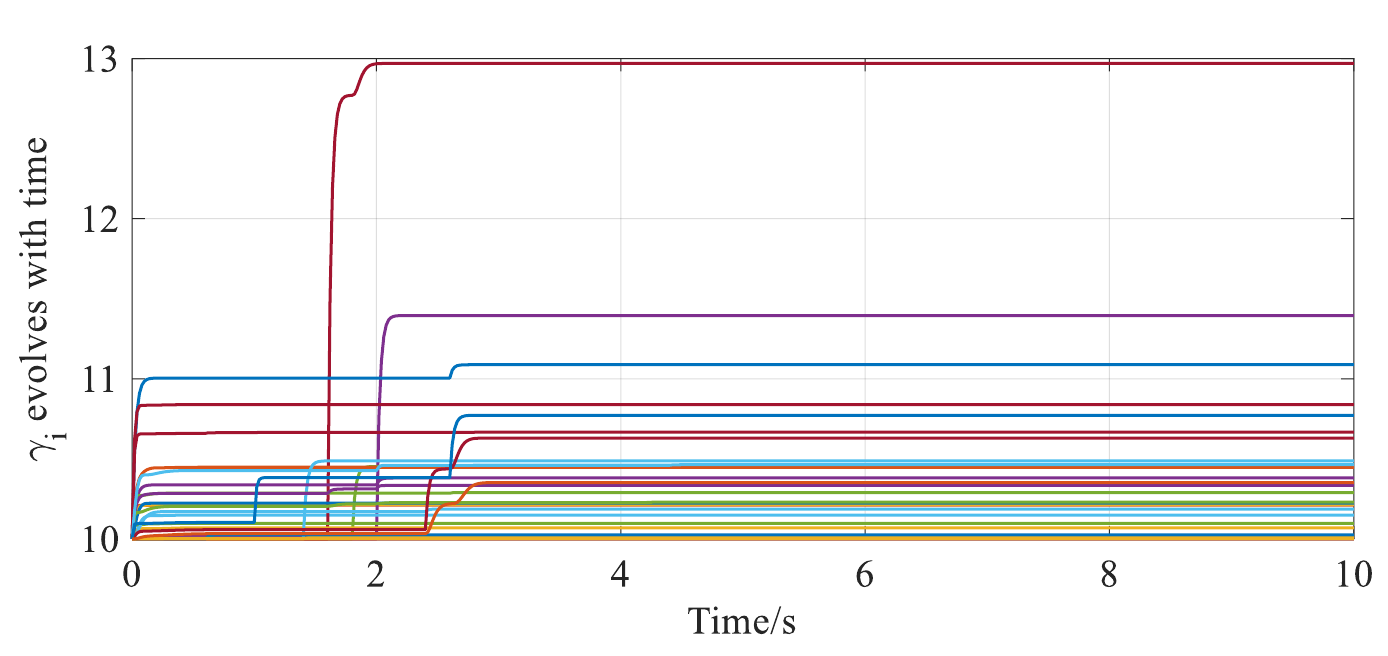}\\
	\caption{Time-varying trajectory of coupling gain}\label{fig-gamma}
\end{figure*}

%
%
%
%

\section{Simulation with a power system}\label{sec6}

To verify the effectiveness of the proposed method, a power system consisting of four power-generation areas is considered. The dynamics of the $i$th area are governed by
\begin{align*}
	&\dot{\chi}_i=A_{ii}\chi_i+\sum_{j=1,~j\neq i}^4A_{ij}\chi_j,\\
	&y_i=C_i\chi_i,
\end{align*}
where $\chi_i\in\mathbb{R}^4$ and $y_i\in\mathbb{R}^2$ are system states and measurement outputs, respectively; $A_{ij}\in\mathbb{R}^{4\times 4}$ for $i,j=1,2,3,4$ are system matrices; and $C_i\in\mathbb{R}^{2\times 4}$ is the output matrix of the $i$th area. $\chi_i$ takes the form $\chi_i=col\{\Delta\theta_i,\Delta w_i,\Delta P_{m_i}, \Delta P_{v_i}\}$. Herein, $\Delta\theta_i$ is the angular displacement of the rotor relative to the stationary reference axis on the stator; $\Delta w_i$ is the change of rotation mass; $\Delta P_{m_i}$ is the change of mechanical power; and $\Delta P_{v_i}$ is the change of steam value position. 

Now, denote $M_i$ the inertia constant, $R_i$ the speed regulation rate, $D_i$ the percentage of load change, $T_{t_i}$ the time constant of the prime mover, $T_{g_i}$ the time constant of the governor, and $P_{ij}$ the slope of the power angle curve between power-generation areas $i$ and $j$ at the initial operating angle.

Then, $A_{ii}$, $A_{ij}$, $C_i$ are expressed as
\begin{align*}
	&A_{ii}=\begin{bmatrix}
		0&1&0&0\\
		-\frac{\sum_{j=1}^NP_{ij}}{2M_i}&-\frac{D_i}{2M_i}&\frac{1}{2M_i}&0\\
		0&0&-\frac{1}{T_{t_i}}&0\\
		0&-\frac{1}{R_iT_{g_i}}&0&-\frac{1}{T_{g_i}}
	\end{bmatrix},\\
	&A_{ij}=\begin{bmatrix}
		0&0&0&0&\\
		\frac{P_{ij}}{2M_i}&0&0&0\\
		0&0&0&0&\\
		0&0&0&0&
	\end{bmatrix},~C_i^T=\begin{bmatrix}
		1&0\\0&0\\0&0\\0&1
	\end{bmatrix}.
\end{align*}

\begin{figure*}[!t]
	\centering
	\includegraphics[width=16cm]{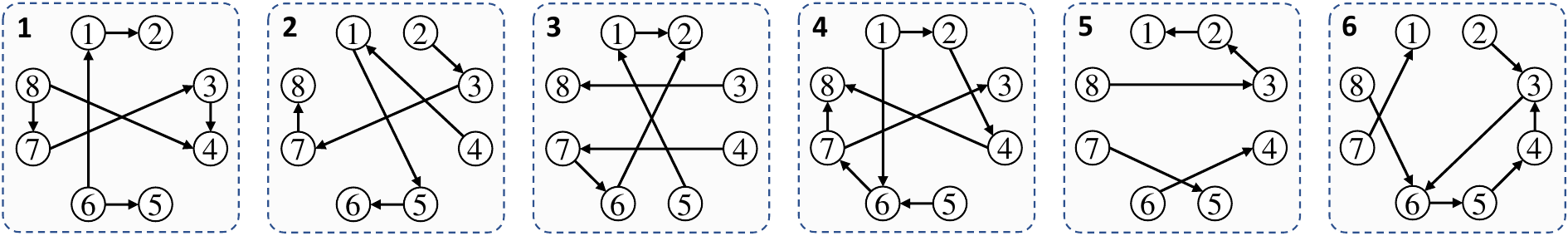}\\
	\caption{Augmented communication networks with $8$ agents}\label{topology2}
\end{figure*}

\begin{table*}\caption{Parameters of the power system}\label{table2}
	\centering
	\renewcommand\arraystretch{1.2}
	\begin{tabular}{p{1cm}<{\centering}p{1.3cm}<{\centering}p{1.3cm}<{\centering}p{1.3cm}<{\centering}p{1.3cm}<{\centering}p{1.3cm}<{\centering}p{1.3cm}<{\centering}p{1.3cm}<{\centering}p{1.3cm}<{\centering}}
		\hline
		~     &Area 1   & Area 2  & Area 3 & Area 4 & Area 5   & Area 6  & Area 7 & Area 8\\
		\hline
		$M_i$    & $12$    & $10$    & $8$    & $8$ & $12$    & $10$    & $8$    & $8$ \\
		$R_i$    & $0.05$  & $0.0625$& $0.08$ & $0.08$ & $0.05$  & $0.0625$& $0.08$ & $0.08$\\
		$D_i$    & $0.7$   & $0.9$   & $0.9$  & $0.7$ & $0.7$   & $0.9$   & $0.9$  & $0.7$ \\
		$T_{t_i}$& $0.65$  & $0.4$   & $0.3$  & $0.6$ & $0.65$  & $0.4$   & $0.3$  & $0.6$ \\
		$T_{g_i}$& $0.1$   & $0.1$   & $0.1$  & $0.1$ & $0.1$   & $0.1$   & $0.1$  & $0.1$ \\
		\hline    
	\end{tabular}
\end{table*}

\begin{figure*}[!t]
	\centering
	\includegraphics[width=18cm]{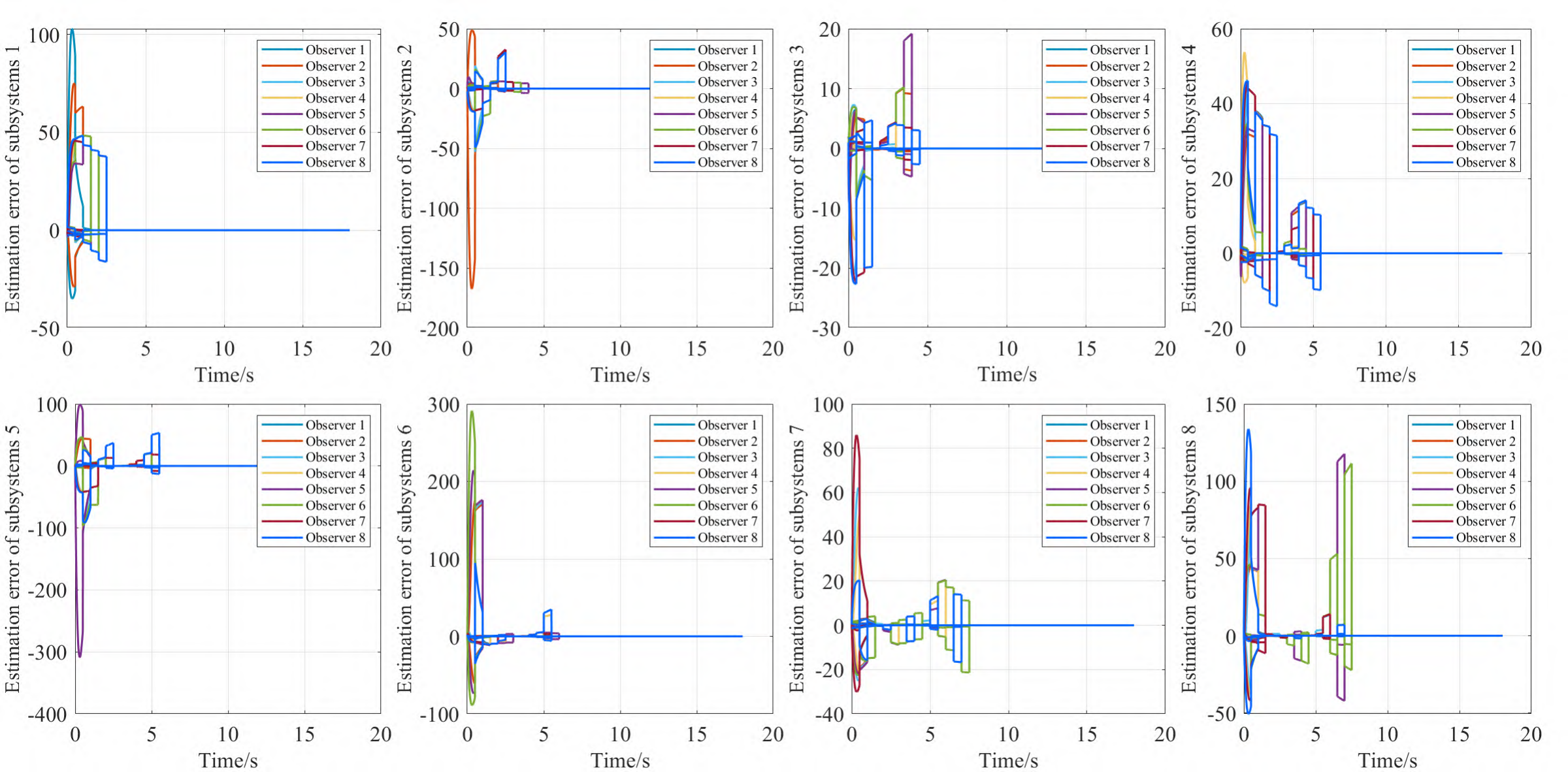}\\
	\caption{Error dynamics of distributed observer with network transformation mapping}\label{fig-powererror8}
\end{figure*}

\begin{figure*}[!t]
	\centering
	\includegraphics[width=18cm]{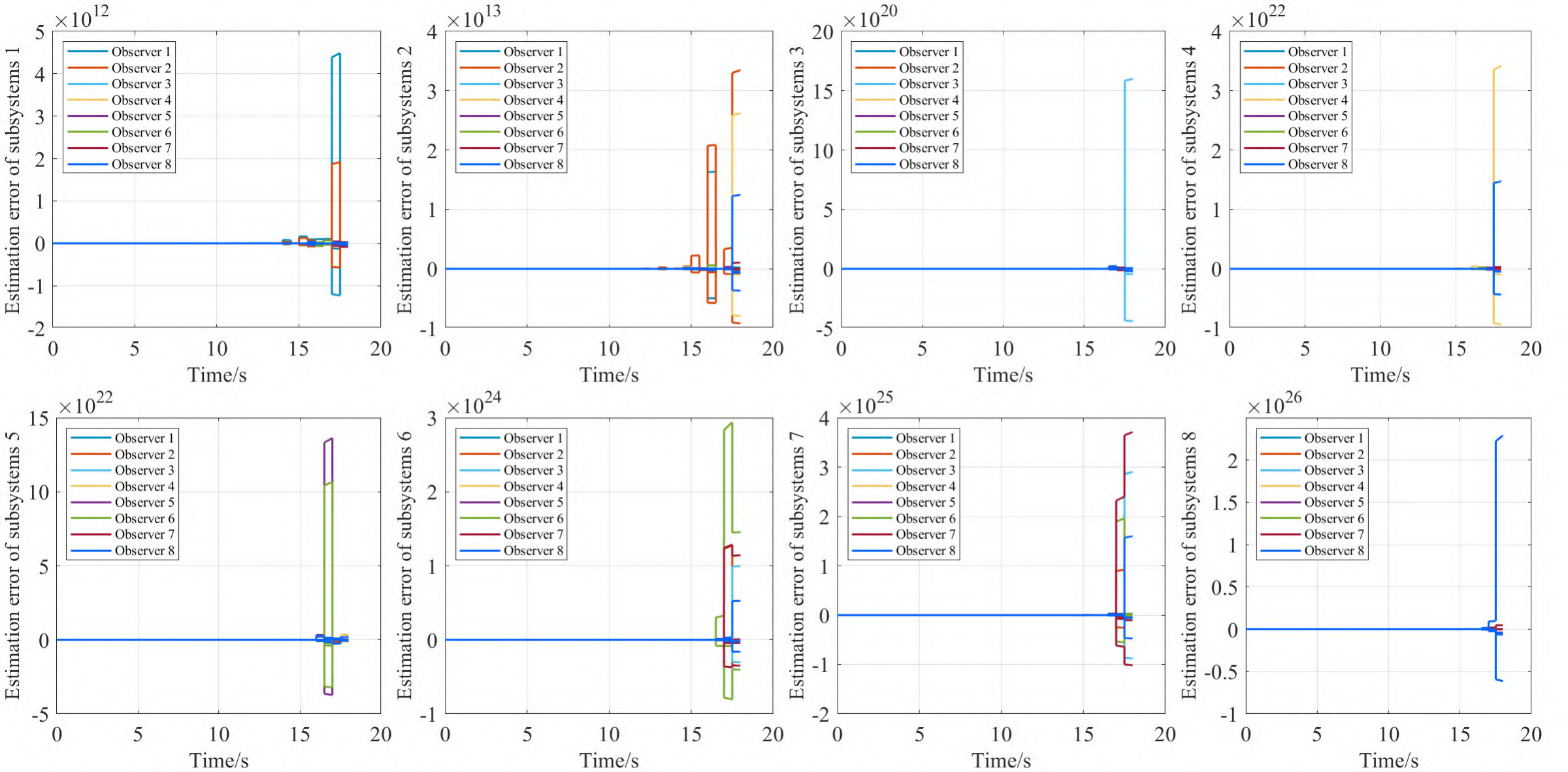}\\
	\caption{Error dynamics of distributed observer without network transformation mapping}\label{fig-powerwithout8}
\end{figure*}

Specific values of the aforementioned parameters are shown in Table \ref{table}. It is assumed that four nodes locate in four areas, and each node has access to $y_i=C_i\chi$ as well as output matrix $C$. The communication networks subject to Assumption \ref{assume3} and \ref{assume2} among the four nodes are shown in Figure \ref{topology}. {\color{red} Note that the eigenvalues of $A=[A_{ij}]_{i,j=1}^N$ are
	\begin{align*}
		&\left\{\underbrace{0,\ldots,0}_{4},\underbrace{-0.029,\ldots,-0.029}_{4},\right.\\
		&\quad\quad\left.\underbrace{-1.54,\ldots,-1.54}_{4},\underbrace{-10,\ldots,-10}_{4}\right\}.
\end{align*}
Therefore, the considered open-loop system is polynominal unstable because it contains multiple zero eigenvalues. This means that the results, which require the condition of open loop stability \cite{zhang2023distributed} and the conditions of undirected networks \cite{Xu2020IFAC,YANG2023110690}, are not well suited for this simulation example.} 

{\color{red}Note that the distributed observer in this example over undirected switching networks can achieve omniscience asymptotically with the approach proposed in \cite{Xu2021TCYB}. However, the method in \cite{Xu2021TCYB} is no longer available when all topologies in the augmented networks are replaced by directed networks. }This paper designs a distributed observer by (\ref{do1}) and (\ref{do2}). The initial values of $\gamma_{i,k}$ in (\ref{doubler}) are chosen to be $10$ as well as $\gamma=100$. $\chi_i(0)$ are set as $[1,1,1,1]^T$ for all $i=1,2,3,4$. The initial values of local observers are chosen randomly in the interval $[-3,3]$. Observer gain $H_{io}$ is chosen so that the eigenvalues of $A_{io}-H_{io}C_{io}$ locate at $\{-12,-13,-14,-15\}$.  

To achieve asymptotic omniscience, we can design distributed observer by (\ref{do1}) and (\ref{do2}) with network transformation mapping. The trajectories of this distributed observer with the fix or adaptive coupling gain are shown in Figure \ref{fig-powererror}, which shows that the distributed observer can achieve omniscience asymptotically under the jointly connected switching directed networks. {\color{red}Comparison results can be seen in Figure \ref{fig-powerwithout}, which shows that the error dynamics cannot stay at zero without using this paper's method (there is no network transformation mapping when designing a distributed observer). This indicates that the current method, including \cite{zhang2023distributed,Xu2020IFAC,YANG2023110690}, are invalid for this example. Moreover, Figure \ref{fig-powerwithout} also implies that Figure \ref{fig-powererror} is not achieved by a fast switching method because Liu. K et al have proved that the distributed observer with jointly connected switching topologies can achieve omniscience asymptotically if the networks switch fast enough \cite{Liu2018Cooperative,2017CooperativeTCS}. Therefore, the instability shown in Figure \ref{fig-powerwithout} indicates that the dwell time ($0.4$s) of each network is longer than the upper limit required for fast switching.} 

Figure \ref{fig-gamma} shows the adaptive coupling gain. {\color{blue}Note that the adaptive coupling gain in this figure converges to a large value. This is because the error dynamics under switching topology decrease to zero with fluctuation. Subsequently, the large error during the fluctuation process will lead to a large value of $\omega_{i,k}$ in equation (\ref{ada-gamma2}). As a result, the adaptive coupling gain increases to a huge value in a short period of time. Nonetheless, since there are no issues in the observer, such as actuator saturation, that pose a threat to system safety, the large coupling gain will not affect system operation. Therefore, although the adaptive strategy may lead to excessive coupling gain, we still adopt this strategy to avoid solving complex inequalities (\ref{gammai}) and (\ref{W_i}).}

{\color{red}To further demonstrate the effectiveness of our theory, we expand the number of subsystems and nodes to $8$. The communication networks between them are shown in Figure \ref{topology2}, and the parameters of the power system containing $8$ areas are given in Table \ref{table2}. We can also calculate that the eigenvalues of $A=[A_{ij}]_{i,j=1}^N$ are
\begin{align*}
	&\left\{\underbrace{0,\ldots,0}_{8},\underbrace{-0.029,\ldots,-0.029}_{8},\right.\\
	&\quad\quad\left.\underbrace{-1.54,\ldots,-1.54}_{8},\underbrace{-10,\ldots,-10}_{8}\right\}.
\end{align*}
Note that the system is also polynomial unstable, so only the distributed observer developed in this paper can achieve omniscience asymptotically for the underlying system. Figure \ref{fig-powererror8} shows the error dynamics of our distributed observer. It indicates that that the error dynamics of the distributed observer fluctuate and converge to $0$ gradually when network transformation mapping exists. In addition, the unstable results shown in Figure \ref{fig-powerwithout8} also demonstrate the indispensability of network transformation mapping.
}

\section{Conclusions}\label{sec7}
In this work, we have concerned with distributed observer under directed switching topologies. To address the challenge whereby there are no strong connected branches in the augmented communication networks, a network transformation mapping has been creatively proposed to construct the spanning tree so that the asymptotic omniscience of the distributed observer can be achieved by a state decomposition and reorganization method. The distributed observer with adaptive coupling gain has been further developed for saving the tedious calculation required by the verification of sufficient conditions. Simulation results have shown that both the fixed coupling gain and the adaptive coupling gain can help the distributed observer achieve omniscience asymptotically under the directed switching topologies. 


\section*{Appendix}
\subsection*{A. Proof of Lemma \ref{cyb}}

According to linear system theory, the solution of a linear time-varying system (\ref{eq-lem1}) is
\begin{align*}
	x(t)=\Phi(t,t_0)x(t_0)+\int_{t_0}^t\Phi(t,\tau)M\xi(\tau){\rm d}\tau,
\end{align*}
where $\Phi(t,t_0)$ is a state transition matrix of $\dot{x}(t)=A(t)x(t)$. Since $\mathbbm{x}_0\leq a_1e^{-a_2(t-t_0)}$, we have $\|\Phi(t,t_0)\|\leq a_1e^{-a_2(t-t_0)}$. Also, $\lim_{t\to\infty}\xi(t)=0$ indicates that there exists $T_1>0$ such that $\|\xi(t)\|\leq \varepsilon(a_2/a_1-\|M\|)/(\|M\|c)\triangleq\varepsilon'$ with $c=\sup_{t>t_0}\|\xi(t)\|$ $\forall t>T_1$ and $\forall\varepsilon>0$. Furthermore, there exists $T_2>0$ such that $e^{-a_2(t-t_0)}\leq \varepsilon'$ $\forall t>T_2$. Then, by defining $x'(t)=\int_{t_0}^t\Phi(t,\tau)M\xi(\tau){\rm d}\tau$, we have
\begin{align*}
	\left\|x'(T)\right\|\leq&\int_{t_0}^{T_m}\|\Phi(t,\tau)\|\|M\|\|\xi(\tau)\|{\rm d}\tau\\
	&+\int_{T_m}^T\|\Phi(t,\tau)\|\|M\|\|\xi(\tau)\|{\rm d}\tau\\
	\leq&\int_{t_0}^{T_m}a_1e^{-a_2(t-\tau)}\|M\|c{\rm d}\tau\\
	&+\int_{T_m}^Ta_1e^{-a_2(t-\tau)}\|M\|\varepsilon'{\rm d}\tau\\
	\leq&\frac{a_1\|M\|c\varepsilon'}{a_2}+\frac{a_1\|M\|\varepsilon'}{a_2}=\varepsilon
\end{align*}
for arbitrary $T>T_m=\max\{T_1,T_2\}$. This implies $\lim_{t\to\infty}x'(t)=0$. Thus, $\lim_{t\to\infty}x(t)=0$.

\subsection*{B. Proof of the convergence with adaptive gain}

In this appendix, we only prove that $\varepsilon_{11}$ and $\xi_{1\star}$ can still converge under the condition of adaptive growth of coupling gain (a similar proof is well suited for $\varepsilon_{i1}$ and $\xi_{i\star}$). However, since the coupling gain is time-varying, the convergence rate is also time-varying, so estimation of the convergence rate that corresponds to different coupling gain values still needs to follow the method in Theorem 1. 

We show the proof as follows. First, we consider the performance of $\varepsilon_{11}$. Since its dynamics
\begin{align*}
	\dot{\varepsilon}_{11}=\left(I_{\pi_1(\sigma)}\otimes (A_{1o}-H_{1o}C_{1o})\right)\varepsilon_{11}
\end{align*}
are stable, a symmetric positive definite matrix $P_{1o}$ can be found such that
\begin{align*}
	P_{1o}(A_{1o}-H_{1o}C_{1o})+(A_{1o}-H_{1o}C_{1o})^TP_{1o}=-2\gamma_1^*I_{v_1}.
\end{align*}
Then, by setting $V_{1,\varepsilon}^a=\frac{1}{2}\varepsilon_{11}^T(I_{\pi_1(\sigma)}\otimes P_{1o})\varepsilon_{11}$, we have $\dot{V}_{1,\varepsilon}^a\leq -\gamma_1^*\varepsilon_{11}^T\varepsilon_{11}$. To move on, the Lyapunov function with respect to $\xi$ can be chosen as 
\begin{align*}
	V_{1,\xi}^a=&\frac{1}{2}\sum_{j\in\mathcal{V}_{1}^{\sigma(t)}\backslash\{1\}}q_{1o,j}(2\gamma_{j,1}+\omega_{j,1})\omega_{j,1}\notag\\
	&+\sum_{j\in\mathcal{V}_{1}^{\sigma(t)}\backslash\{1\}}q_{1o,j}(\gamma_{j,1}-\gamma_1^*)^2.
\end{align*}
It follows that
\begin{align*}
	\dot{V}_{1,\xi}^a=&\sum_{j\in\mathcal{V}_{1}^{\sigma(t)}\backslash\{1\}}q_{1o,j}^\sigma(\gamma_{j,1}+q_{1o,j}^\sigma\omega_{j,1})\dot{\omega}_{j,1}\notag\\
	&+\sum_{j\in\mathcal{V}_{1}^{\sigma(t)}\backslash\{1\}}q_{1o,j}^\sigma(\gamma_{j,1}+\omega_{j,1}-\gamma_1^*)\dot{\gamma}_{j,1}\notag\\
	=&\xi_{1\star}^T\left(\Pi_1^{\sigma}Q_{1o}^\sigma\otimes I_{v_1}\right)\dot{\xi}_{1\star}+\xi_{1\star}^T\left(\Pi_1^{\sigma}Q_{1o}^{\sigma}\otimes I_{v_1}\right)\xi_{1\star}\notag\\
	&-\gamma_1^*\xi_{1\star}^T\left(Q_{1o}^{\sigma}\otimes I_{v_1}\right)\xi_{1\star}\notag\\
	=&\xi_{1\star}^T\left(\Pi_1^{\sigma}Q_{1o}^{\sigma}\otimes A_{1o}\right)\xi_{1\star}\notag\\
	&-\xi_{1\star}^T\left(\Pi_1^{\sigma}Q_{1o}^{\sigma}\mathcal{H}_{\pi_1(\sigma)}\Pi_1^{\sigma}\otimes I_{v_1}\right)\xi_{1\star}\notag\\
	&+\xi_{1\star}^T\left(\Pi_1^{\sigma}Q_{1o}^{\sigma}B_{1\star}^{\sigma(t)}\otimes H_{1o}C_{1o}\right)\varepsilon_{11}\notag\\
	&+\xi_{1\star}^T\left(\Pi_1^{\sigma}Q_{1o}^{\sigma}\otimes I_{v_1}\right)\xi_{1\star}-\gamma_1^*\xi_{1\star}^T\left(Q_{1o}^{\sigma}\otimes I_{v_1}\right)\xi_{1\star}\notag\\
	\leq& \bar{\lambda}(Q_{1o}^{\sigma}\otimes A_{1o})\xi_{1\star}^T\left(\Pi_1^{\sigma}\otimes I_{v_1}\right)\xi_{1\star}\notag\\
	&+\bar{\lambda}(Q_{1o}^{\sigma})\xi_{1\star}^T\left(\Pi_1^{\sigma}\otimes I_{v_1}\right)\xi_{1\star}\notag\\
	&-\gamma_1^*\xi_{1\star}^T\left(\Pi_1^{\sigma}\Pi_1^{\sigma}\otimes I_{v_1}Q_{1o}^{\sigma}\otimes I_{v_1}\right)\xi_{1\star}\notag\\
	&+\bar{\lambda}\left(Q_{1o}^{\sigma}B_{1\star}^{\sigma(t)}\otimes H_{1o}C_{1o}\right)\xi_{1\star}^T\left(\Pi_1^{\sigma}\otimes I_{v_1}\right)\varepsilon_{11}.\notag
\end{align*}
Now, knowledge of mean value inequality leads to
\begin{align*}
	&-\gamma_1^*\xi_{1\star}^T\left(\Pi_1^{\sigma}\Pi_1^{\sigma}\otimes I_{v_1}\right)\xi_{1\star}-\gamma_1^*\xi_{1\star}^T\left(Q_{1o}^{\sigma}\otimes I_{v_1}\right)\xi_{1\star}\notag\\
	=&-\gamma_1^*\sum_{j\in\mathcal{V}_{1}^{\sigma(t)}\backslash\{1\}}(\gamma_{j,1}+\omega_{j,1})^2\xi_{j,1}^T\xi_{j,1}-q_{1o,j}\xi_{j,1}^T\xi_{j,1}\notag\\
	\leq&-2\gamma_1^*\sum_{j\in\mathcal{V}_{1}^{\sigma(t)}\backslash\{1\}}\sqrt{q_{1o,j}}(\gamma_{j,1}+\omega_{j,1})\xi_{j,1}^T\xi_{j,1}\notag\\
	\leq&-2\underline{\lambda}(Q_{1o}^{\sigma})\gamma_1^*\xi_{1\star}^T\left(\Pi_1^{\sigma}\otimes I_{v_1}\right)\xi_{1\star}.
\end{align*}

Let $V_{1,1}^a=V_{1,\varepsilon}^a+V_{1,\xi}^a$, $\wp_{1,1}=-2\underline{\lambda}(Q_{1o}^{\sigma})\gamma_1^*+\bar{\lambda}(Q_{1o}^{\sigma}\otimes A_{1o})+\bar{\lambda}(Q_{1o}^{\sigma})$, $\wp_{1,2}=\bar{\lambda}\left(Q_{1o}^{\sigma}B_{1\star}^{\sigma(t)}\otimes H_{1o}C_{1o}\right)$ and yield
\begin{align*}
	\dot{V}_{1,1}^a\leq& \wp_{1,1}\xi_{1\star}^T\bar{\Pi}_1^\sigma\xi_{1\star}+\wp_{1,2}\xi_{1\star}^T\left(\Pi_1^{\sigma}\otimes I_{v_1}\right)\varepsilon_{11}-\gamma_1^*\varepsilon_{11}^T\varepsilon_{11}\notag\\
	=&\begin{bmatrix}\varepsilon_{11}^T&\xi_{1\star}^T\end{bmatrix}\begin{bmatrix}\wp_{1,1}\bar{\Pi}_1^\sigma&\frac{1}{2}\wp_{1,2}\bar{\Pi}_1^\sigma\\ \frac{1}{2}\wp_{1,2}\bar{\Pi}_1^\sigma&-\gamma_1^* I_{\pi_{1}(\sigma)\times v_1}\end{bmatrix}\begin{bmatrix}\varepsilon_{11}\\ \xi_{1\star}\end{bmatrix},
\end{align*}
where $\bar{\Pi}_1^\sigma=\Pi_1^{\sigma}\otimes I_{v_1}$. In light of the Schur complement lemma, $\dot{V}_{1,1}^a<0$ if and only if 
\begin{align*}
	\wp_{1,1}\left(\Pi_1^{\sigma}\otimes I_{v_1}\right)+\frac{1}{4\gamma_1^*}\wp_{1,2}^2\left(\Pi_1^{\sigma}\otimes I_{v_1}\right)^2<0.
\end{align*}
Since $\Pi_1^{\sigma}\otimes I_{v_1}>0$, there exists a large enough $\gamma_1^*$ such that
\begin{align}\label{ada-gamma2}
	\wp_{1,1}+\frac{1}{4\gamma_1^*}\wp_{1,2}^2\left(\Pi_1^{\sigma}\otimes I_{v_1}\right)<0.
\end{align}
Therefore, $\varepsilon_{11}$ and $\xi_{1\star}$ have a convergence trend because (\ref{ada-gamma2}) leads to $\dot{V}_{1,1}<0$. Based on similar approaches, we can also show that $\varepsilon_{ii}$ and $\xi_{i\star}$ have convergence trend.

\ifCLASSOPTIONcaptionsoff
  \newpage
\fi



%


\bibliographystyle{IEEEtran}
\bibliography{ref-tac}

\end{document}